\documentclass[11pt, a4paper]{article}
\usepackage[cp1251]{inputenc}
\usepackage{amsmath} \usepackage{euscript}
\usepackage{xcolor}
\usepackage{mymatrx6}
\usepackage{graphicx}
\usepackage{longtable}
\usepackage{mathrsfs}
\usepackage{amssymb, amscd}
\def\mymatrix{\MyMatrixwithdelims..}
\Linestrue\Autonumtrue
\usepackage{mathrsfs}
\def\bluecommand{\textcolor[cmyk]{1,0,0,0}}
\def\redcommand{\textcolor[cmyk]{0,1,1,0}}
\def\greencommand{\textcolor[cmyk]{1,0,1,0}}
\def\yellowcommand{\textcolor[cmyk]{0,0.25,1,0}}
\def\graycommand{\textcolor[cmyk]{0,0,0,0.2}}

\begin{document}
	
	\newcounter{bnomer} \newcounter{snomer}
	\newcounter{bsnomer}
	\setcounter{bnomer}{0}
	\renewcommand{\thesnomer}{\thebnomer.\arabic{snomer}}
	\renewcommand{\thebsnomer}{\thebnomer.\arabic{bsnomer}}
	\renewcommand{\refname}{\begin{center}\large{\textbf{References}}\end{center}}
	
	\setcounter{MaxMatrixCols}{14}
	
	\newcommand\restr[2]{{
			\left.\kern-\nulldelimiterspace 
			#1 
			\right|_{#2} 
	}}
	
	\newcommand{\sect}[1]{%
		\setcounter{snomer}{0}\setcounter{bsnomer}{0}
		\refstepcounter{bnomer}
		\par\bigskip\begin{center}\large{\textbf{\arabic{bnomer}. {#1}}}\end{center}}
	\newcommand{\sst}[1]{%
		\refstepcounter{bsnomer}
		\par\bigskip\textsc{\arabic{bnomer}.\arabic{bsnomer}. {#1}}.}
	\newcommand{\defi}[1]{%
		\refstepcounter{snomer}
		\par\medskip\textbf{Definition \arabic{bnomer}.\arabic{snomer}. }{#1}\par\medskip}
	\newcommand{\theo}[2]{%
		\refstepcounter{snomer}
		\par\textbf{Theorem \arabic{bnomer}.\arabic{snomer}. }{#2} {\emph{#1}}\hspace{\fill}$\square$\par}
	\newcommand{\mtheop}[2]{%
		\refstepcounter{snomer}
		\par\textbf{Theorem \arabic{bnomer}.\arabic{snomer}. }{\emph{#1}}
		\par\textsc{Proof}. {#2}\hspace{\fill}$\square$\par}
	\newcommand{\mcorop}[2]{%
		\refstepcounter{snomer}
		\par\textbf{Corollary \arabic{bnomer}.\arabic{snomer}. }{\emph{#1}}
		\par\textsc{Proof}. {#2}\hspace{\fill}$\square$\par}
	\newcommand{\mtheo}[1]{%
		\refstepcounter{snomer}
		\par\medskip\textbf{Theorem \arabic{bnomer}.\arabic{snomer}. }{\emph{#1}}\par\medskip}
	\newcommand{\theobn}[1]{%
		\par\medskip\textbf{Theorem. }{\emph{#1}}\par\medskip}
	\newcommand{\theoc}[2]{%
		\refstepcounter{snomer}
		\par\medskip\textbf{Theorem \arabic{bnomer}.\arabic{snomer}. }{#1} {\emph{#2}}\par\medskip}
	\newcommand{\mlemm}[1]{%
		\refstepcounter{snomer}
		\par\medskip\textbf{Lemma \arabic{bnomer}.\arabic{snomer}. }{\emph{#1}}\par\medskip}
	\newcommand{\mprop}[1]{%
		\refstepcounter{snomer}
		\par\medskip\textbf{Proposition \arabic{bnomer}.\arabic{snomer}. }{\emph{#1}}\par\medskip}
	\newcommand{\theobp}[2]{%
		\refstepcounter{snomer}
		\par\textbf{Theorem \arabic{bnomer}.\arabic{snomer}. }{#2} {\emph{#1}}\par}
	\newcommand{\theop}[2]{%
		\refstepcounter{snomer}
		\par\textbf{Theorem \arabic{bnomer}.\arabic{snomer}. }{\emph{#1}}
		\par\textsc{Proof}. {#2}\hspace{\fill}$\square$\par}
	\newcommand{\theosp}[2]{%
		\refstepcounter{snomer}
		\par\textbf{Theorem \arabic{bnomer}.\arabic{snomer}. }{\emph{#1}}
		\par\textsc{Sketch of the proof}. {#2}\hspace{\fill}$\square$\par}
	\newcommand{\exam}[1]{%
		\refstepcounter{snomer}
		\par\medskip\textbf{Example \arabic{bnomer}.\arabic{snomer}. }{#1}\par\medskip}
	\newcommand{\deno}[1]{%
		\refstepcounter{snomer}
		\par\textbf{Notation \arabic{bnomer}.\arabic{snomer}. }{#1}\par}
	\newcommand{\lemm}[1]{%
		\refstepcounter{snomer}
		\par\textbf{Lemma \arabic{bnomer}.\arabic{snomer}. }{\emph{#1}}\hspace{\fill}$\square$\par}
	\newcommand{\lemmp}[2]{%
		\refstepcounter{snomer}
		\par\medskip\textbf{Lemma \arabic{bnomer}.\arabic{snomer}. }{\emph{#1}}
		\par\textsc{Proof}. {#2}\hspace{\fill}$\square$\par\medskip}
	\newcommand{\coro}[1]{%
		\refstepcounter{snomer}
		\par\textbf{Corollary \arabic{bnomer}.\arabic{snomer}. }{\emph{#1}}\hspace{\fill}$\square$\par}
	\newcommand{\mcoro}[1]{%
		\refstepcounter{snomer}
		\par\textbf{Corollary \arabic{bnomer}.\arabic{snomer}. }{\emph{#1}}\par\medskip}
	\newcommand{\corop}[2]{%
		\refstepcounter{snomer}
		\par\textbf{Corollary \arabic{bnomer}.\arabic{snomer}. }{\emph{#1}}
		\par\textsc{Proof}. {#2}\hspace{\fill}$\square$\par}
	\newcommand{\nota}[1]{%
		\refstepcounter{snomer}
		\par\medskip\textbf{Remark \arabic{bnomer}.\arabic{snomer}. }{#1}\par\medskip}
	\newcommand{\propp}[2]{%
		\refstepcounter{snomer}
		\par\medskip\textbf{Proposition \arabic{bnomer}.\arabic{snomer}. }{\emph{#1}}
		\par\textsc{Proof}. {#2}\hspace{\fill}$\square$\par\medskip}
	\newcommand{\hypo}[1]{%
		\refstepcounter{snomer}
		\par\medskip\textbf{Conjecture \arabic{bnomer}.\arabic{snomer}. }{\emph{#1}}\par\medskip}
	\newcommand{\prop}[1]{%
		\refstepcounter{snomer}
		\par\textbf{Proposition \arabic{bnomer}.\arabic{snomer}. }{\emph{#1}}\hspace{\fill}$\square$\par}
	
	\newcommand{\proof}[2]{%
		\par\medskip\textsc{Proof{#1}}. \hspace{-0.2cm}{#2}\hspace{\fill}$\square$\par\medskip}
	
	\makeatletter
	\def\iddots{\mathinner{\mkern1mu\raise\p@
			\vbox{\kern7\p@\hbox{.}}\mkern2mu
			\raise4\p@\hbox{.}\mkern2mu\raise7\p@\hbox{.}\mkern1mu}}
	\makeatother
	
	\newcommand{\okr}[2]{%
		\refstepcounter{snomer}
		\par\medskip\textbf{{#1} \arabic{bnomer}.\arabic{snomer}. }{\emph{#2}}\par\medskip}

	\newcommand{\Ind}[3]{%
		\mathrm{Ind}_{#1}^{#2}{#3}}
	\newcommand{\Res}[3]{%
		\mathrm{Res}_{#1}^{#2}{#3}}
	\newcommand{\epsi}{\varepsilon}
	\newcommand{\tri}{\triangleleft}
	\newcommand{\Supp}[1]{%
		\mathrm{Supp}(#1)}
	
	\newcommand{\lee}{\leqslant}
	\newcommand{\gee}{\geqslant}
	\newcommand{\reg}{\mathrm{reg}}
	\newcommand{\Ann}{\mathrm{Ann}\,}
	\newcommand{\Cent}[1]{\mathbin\mathrm{Cent}({#1})}
	\newcommand{\PCent}[1]{\mathbin\mathrm{PCent}({#1})}
	\newcommand{\Irr}[1]{\mathbin\mathrm{Irr}({#1})}
	\newcommand{\Exp}[1]{\mathbin\mathrm{Exp}({#1})}
	\newcommand{\empr}[2]{[-{#1},{#1}]\times[-{#2},{#2}]}
	\newcommand{\sreg}{\mathrm{sreg}}
	\newcommand{\ilm}{\varinjlim}
	\newcommand{\plm}{\varprojlim}
	\newcommand{\codim}{\mathrm{codim}\,}
	\newcommand{\chara}{\mathrm{char}\,}
	\newcommand{\rk}{\mathrm{rk}\,}
	\newcommand{\chr}{\mathrm{ch}\,}
	\newcommand{\Ker}{\mathrm{Ker}\,}
	\newcommand{\id}{\mathrm{id}}
	\newcommand{\Ad}{\mathrm{Ad}}
	\newcommand{\col}{\mathrm{col}}
	\newcommand{\row}{\mathrm{row}}
	\newcommand{\high}{\mathrm{high}}
	\newcommand{\low}{\mathrm{low}}
	\newcommand{\pho}{\hphantom{\quad}\vphantom{\mid}}
	\newcommand{\fho}[1]{\vphantom{\mid}\setbox0\hbox{00}\hbox to \wd0{\hss\ensuremath{#1}\hss}}
	\newcommand{\wt}{\widetilde}
	\newcommand{\wh}{\widehat}
	\newcommand{\ad}[1]{\mathrm{ad}_{#1}}
	\newcommand{\tr}{\mathrm{tr}\,}
	\newcommand{\GL}{\mathrm{GL}}
	\newcommand{\SL}{\mathrm{SL}}
	\newcommand{\SO}{\mathrm{SO}}
	\newcommand{\Or}{\mathrm{O}}
	\newcommand{\Sp}{\mathrm{Sp}}
	\newcommand{\Sa}{\mathrm{S}}
	\newcommand{\Ua}{\mathrm{U}}
	\newcommand{\Mat}{\mathrm{Mat}}
	\newcommand{\Stab}{\mathrm{Stab}}
	\newcommand{\htt}{\mathfrak{h}}
	\newcommand{\spt}{\mathfrak{sp}}
	\newcommand{\slt}{\mathfrak{sl}}
	\newcommand{\sot}{\mathfrak{so}}

	\newcommand{\vfi}{\varphi}
	\newcommand{\aad}{\mathrm{ad}}
	\newcommand{\vpi}{\varpi}
	\newcommand{\teta}{\vartheta}
	\newcommand{\Bfi}{\Phi}
	\newcommand{\Fp}{\mathbb{F}}
	\newcommand{\Rp}{\mathbb{R}}
	\newcommand{\Zp}{\mathbb{Z}}
	\newcommand{\Cp}{\mathbb{C}}
	\newcommand{\Ap}{\mathbb{A}}
	\newcommand{\Pp}{\mathbb{P}}
	\newcommand{\Kp}{\mathbb{K}}
	\newcommand{\Np}{\mathbb{N}}
	\newcommand{\ut}{\mathfrak{u}}
	\newcommand{\at}{\mathfrak{a}}
	\newcommand{\glt}{\mathfrak{gl}}
	\newcommand{\hei}{\mathfrak{hei}}
	\newcommand{\nt}{\mathfrak{n}}
	\newcommand{\mt}{\mathfrak{m}}
	\newcommand{\rt}{\mathfrak{r}}
	\newcommand{\rad}{\mathfrak{rad}}
	\newcommand{\bt}{\mathfrak{b}}
	\newcommand{\unt}{\underline{\mathfrak{n}}}
	\newcommand{\gt}{\mathfrak{g}}
	\newcommand{\vt}{\mathfrak{v}}
	\newcommand{\pt}{\mathfrak{p}}
	\newcommand{\Xt}{\mathfrak{X}}
	\newcommand{\Po}{\mathcal{P}}
	\newcommand{\PV}{\mathcal{PV}}
	\newcommand{\Uo}{\EuScript{U}}
	\newcommand{\Fo}{\EuScript{F}}
	\newcommand{\Do}{\EuScript{D}}
	\newcommand{\Eo}{\EuScript{E}}
	\newcommand{\Jo}{\EuScript{J}}
	\newcommand{\Iu}{\mathcal{I}}
	\newcommand{\Mo}{\mathcal{M}}
	\newcommand{\Nu}{\mathcal{N}}
	\newcommand{\Ro}{\mathcal{R}}
	\newcommand{\Co}{\mathcal{C}}
	\newcommand{\Ko}{\mathcal{K}}
	\newcommand{\So}{\mathcal{S}}
	\newcommand{\Lo}{\mathcal{L}}
	\newcommand{\Ou}{\mathcal{O}}
	\newcommand{\Uu}{\mathcal{U}}
	\newcommand{\Tu}{\mathcal{T}}
	\newcommand{\Au}{\mathcal{A}}
	\newcommand{\Vu}{\mathcal{V}}
	\newcommand{\Du}{\mathcal{D}}
	\newcommand{\Bu}{\mathcal{B}}
	\newcommand{\Sy}{\mathcal{Z}}
	\newcommand{\Sb}{\mathcal{F}}
	\newcommand{\Gr}{\mathcal{G}}
	\newcommand{\Xu}{\mathcal{X}}
	\newcommand{\rtc}[1]{C_{#1}^{\mathrm{red}}}

	\newcommand{\JSpec}[1]{\mathrm{JSpec}\,{#1}}
	\newcommand{\MSpec}[1]{\mathrm{MSpec}\,{#1}}
	\newcommand{\PSpec}[1]{\mathrm{PSpec}\,{#1}}
	\newcommand{\APbr}[1]{\mathrm{span}\{#1\}}
	\newcommand{\APbre}[1]{\langle #1\rangle}
	\newcommand{\APro}[1]{\setcounter{AP}{#1}\Roman{AP}}\newcommand{\apro}[1]{{\rm\setcounter{AP}{#1}\roman{AP}}}
	\newcommand{\ot}{\xleftarrow[]{}}
	\newcounter{AP}
	

	\author{Mikhail Ignatev\and Mikhail Venchakov}
	\date{}
	\title{Characters of the unitriangular group and the Mackey method}\maketitle
	\begin{abstract} Let $U$ be the unitriangular group over a finite field. We consider an interesting class of irreducible complex characters of $U$, so-called characters of depth 2. This is a next natural step after characters of maximal and submaximal dimension, whose description is already known. We explicitly describe the support of a character of depth 2 by a system of defining algebraic equations. After that, we calculate the value of such a character on an element from the support. The main technical tool used in the proofs is the Mackey little group method for semidirect products.
		
		\medskip\noindent{\bf Keywords:} coadjoint orbit, orbit method, irreducible character, Mackey method, polarization.\\
		{\bf AMS subject classification:} 20C15, 17B08, 20D15.\end{abstract}
	
	\begin{flushright}
		\emph{Dedicated to the memory of N.A. Vavilov}
	\end{flushright}
	\vspace{-0.9cm}
	\sect{Introduction}
	
	The\let\thefootnote\relax\footnote{The work of the first author is an output of a research project implemented as part of the Basic Research Program at the National Research University Higher School of Economics (HSE University). The work of the second author was supported by the Russian Science Foundation under grant no. 22--71--10001.} main tool in representation theory of a unipotent group $U$ over a finite field $\Fp_q$ is the orbit method created in 1962 by A.A. Kirillov \cite{Kirillov62}, \cite{Kirillov04}. According to this method, one should consider the Lie algebra $\ut$ of $U$ and the dual space $\ut^*$. The group $U$ acts on $\ut$ via the adjoint action; the dual action of $U$ on $\ut^*$ is called \emph{coadjoint}. It turns out that the complex finite-dimensional irreducible representations of $U$ are in one-to-one correspondence with the coadjoint orbits on $\ut^*$ \cite{Kazhdan77}, see Section~\ref{sect:orbit_method} for the detail. Furthermore, for an orbit $\Omega\subset\ut^*$, one can (in principle) compute the value of the corresponding irreducible character $\chi=\chi_{\Omega}$ on an arbitrary element of $U$, see formula (\ref{formula:chi_Omega}) below. In particular, the \emph{degree} $\deg\chi$ of $\chi$, i.e., the dimension of the corresponding representation, is equal to~$q^{\dim\Omega/2}$.
	
	Let $U=U_n(\Fp_q)$ be the unitriangular group, i.e., the group of $n\times n$ strictly lower triangular matrices with 1's on the diagonal over $\Fp_q$. A complete description of the coadjoint orbits for the group~$U$ is a wild problem. Hence, an interesting question is how to describe ``the most important'' (in some sense) classes of orbits, f.e., orbits of maximal or submaximal dimension. Moreover, formula~(\ref{formula:chi_Omega}) for the character $\chi$ is not explicit, but requires the summation over all elements of the orbit $\Omega$ instead. Thus, even if a description of an orbit is known, the explicit calculation of the corresponding character is an interesting problem itself. Precisely, for an irreducible character $\chi$ of $U$, its \emph{support} is by definition the set $\Supp{\chi}=\{g\in U\mid\chi(g)\neq0\}.$ Clearly, it is a disjoint union of certain conjugacy classes in $U$. So, the problem is as follows: given a coadjoint orbit $\Omega$, describe explicitly conjugacy classes contained in $\Supp{\chi_{\Omega}}$, and then calculate the values of $\chi_{\Omega}$ on such classes.
	
	Orbits of maximal possible dimension $N=2((n-2)+(n-4)+\ldots)$ were classified in the first Kirillov's paper \cite{Kirillov62} on the orbit method (see also \cite[Theorem 3.1]{IgnatevPanov09} for the detailed proof); we recall this description in Example~\ref{exam:0_1_regular_orbits} (i) below. The corresponding irreducible characters were explicitly computed by C. Andr\`e in \cite{Andre01}. We will refer to these characters as characters \emph{of depth} 0. Orbits of submaximal dimension $N-2$ were classified by A.N. Panov in \cite{IgnatevPanov09}. The corresponding characters were calculated by the first author in \cite{Ignatev09}. We will refer to some of these characters as characters \emph{of depth} 1, see the details in Example~\ref{exam:0_1_regular_orbits} (ii) below.
	
	The goal of the present paper is to show how the well-known Mackey little group method allows to perform the next natural step, i.e., to consider the characters \emph{of depth}~2. To give precise definitions, we need some more notation. The Lie algebra $\ut=\ut_n(\Fp_q)$ is the algebra of strictly lower triangular matrices with zeroes on the diagonal. Using the trace form, one can identify~$\ut^*$ with the space $\ut^t$ of strictly upper triangular matrices with zeroes on the diagonal. There is a natural stratification of $\ut^*$ by the $U$-invariant subspaces $\Xu_i$, $0\leq i\leq n-1$, defined as
	\begin{equation*}
		\Xu_i=\{\lambda\in\ut^*\mid\lambda_{j,n}=0\text{ for }j<i\}.
	\end{equation*}
	In particular, $\Xu_0=\ut^*$, and $\Xu_{n-1}$ is naturally isomorphic to the dual space of $\ut_{n-1}(\Fp_q)$.
	
	\defi{We say that a character $\chi_{\Omega}$ is \emph{of depth} $i$, $1\leq i\leq n-1$, if $\Omega$ is contained in $\Xu_i$ and has maximal possible dimension among all orbits lying in $\Xu_i$. (We call such an orbit $i$-\emph{regular}.)}
	
	For instance, the 0-regular orbits are exactly the orbits of maximal dimension $N$ in $\ut^*$, while\break 1-regular orbits has submaximal dimension $N-2$ (but not all orbits of dimension $N-2$ are\break 1-regular). In general, Panov computed the dimension of an $i$-regular orbit: it equals $N-2i$\break \cite[Proposition 3.2]{IgnatevPanov09}. It is well known that the set of linear forms from $\Xu_i$ whose orbits have maximal possible dimension, is open and dense in $\Xu_i$ \cite[\S2.6]{Kraft00}. So, characters of depth $i$ are in some sense ``in general position'' at the corresponding layer of this stratification. Recall again that characters of depth 0 and 1 were already calculated in precise form.
	
	Our main result, formulated and proved in Sections \ref{sect:main_result} and \ref{sect:proofs_classes} respectively, is an explicit formula for a character $\chi$ of depth~2. Namely, in Theorem~\ref{theo:2_depth_support} we present defining equations for the support of~$\chi$, while in Theorem~\ref{theo:2_depth_value}, given an element $g\in\Supp{\chi}$, the value $\chi(g)$ is calculated.
	
	The structure of the paper is as follows. In the next section, we briefly recall mail definitions and facts from the orbit method, and then present a classification of 0-regular and 1-regular orbits. Section~\ref{sect:Mackey_method} contains a brief description of the Mackey's little group method in representation theory of semidirect products, as well as its application to the calculation of the characters of depths 0 and~1. In Section~\ref{sect:2_orbits}, we classify all 2-regular orbits in the spirit of paper \cite{IgnatevPanov09}. Next, in Section~\ref{sect:main_result}, we formulate our main results about characters of depth 2, which are finally proved in Sections~\ref{sect:proofs_classes} and~\ref{sect:support} by methods from Sections~\ref{sect:orbit_method} and \ref{sect:Mackey_method}.
	
	\medskip\textsc{Acknowledgements}. We thank Prof. A.N. Panov for fruitful discussions and suggestions. We also thank M.Yu. Panov for his useful suggestions about drawing pictures with \TeX.
	
	\sect{The orbit method}\label{sect:orbit_method}
	
	Throughout this section, $\Kp$ is the algebraic closure of the field $\Fp_q$ with $q$ elements of sufficiently large characteristic $p$, $U(\Kp)$ is a connected unipotent algebraic group over $\Kp$ defined over $\Fp_q$, and $U$ is the set of its $\Fp_q$-points. Below we briefly present the orbit method, which describes the irreducible complex representations of $U$ in terms of the coadjoint orbits of this group. This description will be applied for the calculation of a formula for a character of depth 2.
	
	By $\ut(\Kp)$ and $\ut$ we denote the Lie algebras of the groups $U(\Kp)$ and $U$ respectively, and by $\ut(\Kp)^*$ and $\ut^*$ we denote the corresponding dual spaces. Note that one can consider $\ut$ and $\ut^*$ as the sets of $\Fp_q$-points of $\ut(\Kp)$ and $\ut(\Kp)^*$ respectively. The groups $U(\Kp)$, $U$ acts on their Lie algebras via the adjoint action; the dual action on the respective dual spaces is called \emph{coadjoint}. Note that, being an orbit of a connected unipotent group on an affine space, any coadjoint $U(\Kp)$-orbit $\Omega(\Kp)$ on $\ut(\Kp)^*$ is an affine subvariety; in fact, it is isomorphic to an affine space. It is known that its dimension $\dim\Omega(\Kp)$ is even. Moreover, given a linear form $f\in\ut^*\subset\ut(\Kp)^*$, denote by $\Omega_f$, $\Omega_f(\Kp)$ its coadjoint orbits under the action of $U$, $U(\Kp)$ respectively, then $\Omega_f$ is the set of $\Fp_q$-points of $\Omega_f$, so $|\Omega_f|=q^{\dim\Omega_f(\Kp)}$. Throughout the paper, we will call $\dim\Omega_f(\Kp)$ the \emph{dimension} of the orbit $\Omega_f$ and denote it by $\dim\Omega_f$.
	
	According to the orbit method, to each linear form $\lambda\in\ut^*$ one can attach the irreducible representation of $U$. To do this, we need some more definitions. Let $\beta_{\lambda}$ be a bilinear form on $\ut$ defined by $\beta_{\lambda}(x,y)=\lambda([x,y])$ for $x,y\in\ut$.
	
	\defi{A subalgebra $\pt\subset\ut$ is called a \emph{polarization} for $\lambda$ if $\pt$ is a maximal $\beta_{\lambda}$-isotropic subspace. It is well known that $\dim\Omega_{\lambda}=2\codim_{\ut}{\pt}$. According to M. Vergne's result \cite{Vergne}, such a polarization always exists.}
	
	Since the characteristic $p$ is sufficiently large, the usual exponential map $\exp\colon\ut\to U$ is well defined; actually, it is an isomorphism of affine varieties. For example, if $U=U_n(\Fp_q)$, then this map is given by the usual formula $$\exp(x)=\sum_{k=0}^{\infty}x^k/k!=\sum_{k=0}^{n-1}x^k/k!,$$ the latter equality being true because $x^n=0$ for all $x\in\ut$ (here we assume that $p\geq n$). We denote the inverse map by $\ln\colon U\to\ut$. If $g.\lambda$ stands for the result of the coadjoint action for $g\in U$, $\lambda\in\ut^*$, then by definition
	\begin{equation*}
		((\exp x).\lambda)(y)=\lambda\left(\sum_{k=0}^{\infty}\ad{-x}^k(y)/k!\right),~x,y\in\ut.
	\end{equation*}
	Here $\ad{x}\colon\ut\to\ut$, $y\mapsto[x,y]$, is the adjoint operator corresponding to $x$, and the sum in the right-hand side is in fact finite. (Of course, the same formula is true over $\Kp$.)
	
	For example, if $U=U_n(\Fp_q)$, then one can naturally identify $\ut^*$ with the space $\ut^t$ of upper triangular matrices with zeroes on the diagonal by putting $\lambda(x)=\tr(\lambda x)$ for $\lambda\in\ut^t$, $x\in\ut$, as it was mentioned in the Introduction. Under this identification, coadjoint action has the following simple form:
	\begin{equation*}
		g.\lambda=(g\lambda g^{-1})_{\high},
	\end{equation*}
	where $a_{\high}$ is the matrix obtained from an $n\times n$ matrix $a$ by replacing its elements on the diagonal and below by zeroes.
	
	Now, fix a non-trivial homomorphism $\theta\colon\Fp_q\to\Cp^{\times}$, where $\Fp^{\times}$ means the multiplicative group of a field~$\Fp$, and denote $P=\exp\pt$. Since the Baker--Campbell--Hausdorff formula is satisfied, $P$ is a subgroup of~$U$. Furthermore, the condition that $\pt$ is $\beta_{\lambda}$-isotropic implies immediately that $$\psi=\theta\circ\lambda\circ\ln\colon P\to\Cp^{\times}$$ is a one-dimensional complex representation of the group $P$.
	
	We denote by $T$ (respectively, by $\chi$) the induced representation of $U$ (respectively, its character). It turns out that $\chi$ can be computed by the following formula:
	\begin{equation}
		\chi(g)=\dfrac{1}{\sqrt{|\Omega_{\lambda}|}}\sum_{\mu\in\Omega_{\lambda}}\theta(\mu(\ln g)),~g\in U.\label{formula:chi_Omega}
	\end{equation}
	(Note that $\sqrt{|\Omega_{\lambda}|}=q^{\dim\Omega_{\lambda}/2}$.) In particular, $\chi$ depends only on the coadjoint orbit $\Omega=\Omega_{\lambda}$, but not on the specific point $\lambda\in\Omega$ and its polarization $\pt$. We then denote $\chi=\chi_{\Omega}$ and $T=T_{\Omega}$. It follows immediately from formula (\ref{formula:chi_Omega}) that $$\deg\chi_{\Omega}=\dim_{\Cp}T_{\Omega}=\chi_{\Omega}(e)=q^{\dim\Omega_{\lambda}/2}=\sqrt{|\Omega|},$$ where $e$ is the neutral element of $U$.
	
	The orbit method can by shortly expressed as follows.
	\theop{Let $U$\textup, $\ut$\textup, $\ut^*$\textup, etc. be as above. Then \textup{i)} the representation $T_{\Omega}$ is irreducible\textup;\break \textup{ii)} each irreducible representation of $U$ has the form $T_{\Omega}$ for a certain coadjoint orbit $\Omega\subset\ut^*$\textup;\break \textup{iii)} $T_{\Omega_1}$ is isomorphic to $T_{\Omega_2}$ if and only if $\Omega_1=\Omega_2$.}{All these three claims can be easily deduced from formula (\ref{formula:chi_Omega}).}\newpage
	
	Thus, the irreducible representations of the group $U$ bijectively correspond to the coadjoint orbits on $\ut^*$. Unfortunately, even for the group $U=U_n(\Fp_q)$, a complete description of all coadjoint orbits is a wild problem, so one can concentrate on a description of special important classes of orbits.
	
	From now on and to the end of the paper, $U=U_n(\Fp_q)$. It is convenient to introduce the following notation. Let $\Phi(n)$ be the \emph{set of positive roots}, which can be identified with the set of pairs
	\begin{equation*}
		\Phi(n)=\{(i,j),~1\leq j<i\leq n\}.
	\end{equation*}
	Clearly, $\{e_{i,j},~(i,j)\in\Phi(n)\}$ is a basis of $\ut$, where $e_{i,j}$ denotes the usual $(i,j)$th elementary matrix. We denote by $\{e_{i,j}^*,~(i,j)\in\Phi\}$ the dual basis of $\ut^*$; under the identification $\ut^*\cong\ut^t$ defined above, $e_{i,j}^*$ is nothing but $e_{j,i}$. Given a linear form $\lambda\in\ut^*$, we will denote by $\Supp{\lambda}$ its \emph{support}, i.e., the set $$\Supp{\lambda}=\{(i,j)\in\Phi(n)\mid\lambda(e_{i,j})\neq0\},$$
	so that $$\lambda=\sum_{(i,j)\in\Supp{\lambda}}\lambda(e_{i,j})e_{j,i}.$$
	Given a subset $D\subset\Phi(n)$ and a map $\xi\colon D\to\Fp_q^{\times}$, we put $$f_{D,\xi}=\sum_{(i,j)\in D}\xi(i,j)e_{j,i},$$ then $D=\Supp{f_{D,\xi}}$. Finally, to each root $(i,j)\in\Phi(n)$ we assign the set of its \emph{singular roots}: by definition, it is the set $S(i,j)=S^+(i,j)\cup S^-(i,j)$, where $$S^+(i,j)=\{(i,l),~j<l<i\},~S^-(i,j)=\{(l,j),~j<l<i\}.$$
	
	We need some more notation. Consider a partial order on $\Phi(n)$ defined as follows: $(a,b)\geq(c,d)$ if $a\geq c$ and $b\leq d$. For a root $(i,j)\in\Phi(n)$, we call $\row(i,j)=i$ (respectively, $\col(i,j)=j$) the \emph{row} (respectively, the \emph{column}) of this root. Given $l$, we call the sets
	\begin{equation*}
		\Ro_l=\{(l,j),~1\leq j<l\},~\Co_l=\{(i,l),~l<i\leq n\}
	\end{equation*}
	the $l$-th \emph{row}, \emph{column} of $\Phi(n)$ respectively.
	
	Assume that a subset $D\subset\Phi(n)$ satisfies $|D\cap\Ro_l|\leq1$, $|D\cap\Co_l|\leq1$ for all $l$ (such a subset is called a \emph{rook placement}.) For a root $(r,s)\in\Phi(n)\setminus\left(\bigcup_{(i,j)\in D}S(i,j)\right)$, put $$D_{r,s}=\{(i,j)\in D\mid(i,j)>(r,s)\}.$$ Denote $D_{r,s}\cup\{(r,s)\}=\{(i_1,j_1),\ldots,(i_t,j_t)\}$, where $j_1<\ldots<j_t$. For a matrix $\lambda\in\ut^t$, we denote by $\Delta_D^{r,s}(\lambda)$ the minor of the matrix $\lambda$ with the set of rows $j_1,\ldots,j_t$ and the set of columns $\sigma(i_1),\ldots,\sigma(i_t)$, where $\sigma$ is the unique permutation such that $\sigma(i_1)<\ldots<\sigma(i_t)$\label{page:delta}. We are now ready to describe all 0-regular and 1-regular orbits.
	
	\exam{i) Denote by $D^0$ the set \label{exam:0_1_regular_orbits} $\{(n,1),~(n-1,2),~\ldots,~(n-n_0+1,n_0)\}$, where\break $n_0=[(n-1)/2]$. If $n$ is even then this set without the last root is also denoted by $D^0$. Fix also a map $\xi\colon D^0\to\Fp_q^{\times}$, and put $$\Mo=\{(i,j)\in\Phi(n)\mid i>n-j+1\}=\bigcup_{(r,s)\in D^0}S^-(r,s).$$ For instance, let $n=6$. On the picture below we schematically drew this sets: the roots from $D^0$ are marked by $\otimes$, while the roots from $\Mo$ are marked by minuses. We also marked by pluses the roots from $$\bigcup_{(r,s)\in D^0}S^+(r,s)=\{(i,j)\in\Phi(n)\mid i<n-j+1\}=\Phi(n)\setminus(D^0\cup\Mo).$$ Here we identify $\Phi(n)$ with the lower-triangular chessboard of size $n\times n$ in an obvious way: we send a root $(i,j)$ to the $(i,j)$th box.
		\begin{equation*}
			\mymatrix{
				\Bot{2pt}\pho& \pho& \pho& \pho& \pho& \pho\\
				\Rt{2pt}+& \Bot{2pt}\pho& \pho& \pho& \pho& \pho\\
				+& \Rt{2pt}+& \Bot{2pt}\pho& \pho& \pho& \pho\\
				+& +& \Rt{2pt}\otimes& \Bot{2pt}\pho& \pho& \pho\\
				+& \otimes& -& \Rt{2pt}-& \Bot{2pt}\pho& \pho\\
				\otimes& -& -& -& \Rt{2pt}-& \pho\\}
		\end{equation*}
		
		It turned out \cite{Kirillov62} that the orbit $\Omega_{D^0,\xi}=\Omega_{f_{D^0,\xi}}$ is 0-regular, and each 0-regular orbit contains exactly one point of the form $f_{D^0,\xi}$. (Recall that the 0-regular orbits are exactly the orbits of maximal dimension $N$.) Moreover, the subspace $\pt$ of $\ut$ spanned by all vectors $e_{i,j}$, $(i,j)\in\Phi(n)\setminus\Mo$, is a polarization for $f_{D^0,\xi}$. The proof of the latter fact is easy and very similar to the proof of Lemma~\ref{lemm:polarization} below, so we omit it here.
		
		We will also present a system of defining equations for the orbit $\Omega_{D^0,\xi}$. Precisely, $\mu\in\ut^*$ belongs to the orbit $\Omega_{D^0,\xi}$ if and only if $$\Delta_{D^0}^{r,s}(\mu)=\Delta_{D^0}^{r,s}(f_{D^0,\xi})=(-1)^k\prod_{(i,j)\in D^0_{r,s}}\xi(i,j)\text{ for all }(r,s)\in D^0,\text{ where }k=[|D^0_{r,s}|/2].$$
		For instance, in the example above the orbit $\Omega_{D^0,\xi}$ is defined by the following equations:
		\begin{equation*}
			\mu_{1,6}=\xi(6,1),~
			\begin{vmatrix}
				\mu_{1,5}&\mu_{2,5}\\
				\mu_{1,6}&\mu_{2,6}
			\end{vmatrix}=-\xi(6,1)\xi(5,2),~
			\begin{vmatrix}
				\mu_{1,4}&\mu_{2,4}&\mu_{3,4}\\
				\mu_{1,5}&\mu_{2,5}&\mu_{3,5}\\
				\mu_{1,6}&\mu_{2,6}&\mu_{3,6}\\
			\end{vmatrix}=-\xi(6,1)\xi(5,2)\xi(4,3).
		\end{equation*}
		
		ii) To classify all 1-regular coadjoint orbits, denote $$D^1=(D^0\setminus\{(n,1),(n-1,2)\})\cup\{(n-1,1),(n,2)\}\cup R_1,$$ where $R_1$ is either empty or equals $\{(n,n-1)\}$. In the spirit of the previous example, denote $$\Mo=\left(\bigcup_{(r,s)\in D^1}S^-(r,s)\right)\setminus\{(n,n-1)\}.$$
		On the picture below, for $n=8$, we denote the roots from $D^1\setminus R_1$ by $\otimes$, the roots from $\Mo$ by minuses, the root $(n,n-1)$ by square, the root $(n,1)$ by dot, and all other roots by pluses.
		\begin{equation*}
			\mymatrix{
				\Bot{2pt}\pho& \pho& \pho& \pho& \pho& \pho& \pho& \pho\\
				\Rt{2pt}+& \Bot{2pt}\pho& \pho& \pho& \pho& \pho& \pho& \pho\\
				+& \Rt{2pt}+& \Bot{2pt}\pho& \pho& \pho& \pho& \pho& \pho\\
				+& +& \Rt{2pt}+& \Bot{2pt}\pho& \pho& \pho& \pho& \pho\\
				+& +& +& \Rt{2pt}\otimes& \Bot{2pt}\pho& \pho& \pho& \pho\\
				+& +& \otimes& -& \Rt{2pt}-& \Bot{2pt}\pho& \pho& \pho\\
				\otimes& -& -& -& -& \Rt{2pt}-& \Bot{2pt}\pho& \pho\\
				\bullet& \otimes& -& -& -& -& \Rt{2pt}\square& \pho\\}
		\end{equation*}
		
		Fix a map $\xi\colon D^1\to\Fp_q^{\times}$. As it was proved in \cite[Theorem 3.3]{IgnatevPanov09}, the orbit $\Omega_{D^1,\xi}=\Omega_{f_{D^1,\xi}}$ is 1-regular, and each 1-regular orbit contains exactly one point of the form $f_{D^1,\xi}$. (Each 1-regular orbit has submaximal dimension $N-2$, but not all orbits of dimension $N-2$ are 1-regular.) Further, the subspace $\pt\subset\ut$ spanned by all vectors $e_{i,j}$, $(i,j)\in\Phi(n)\setminus\Mo$, is a polarization for $f_{D^1,\xi}$. As in step~(i), the proof of the latter fact is vary similar to the proof of Lemma~\ref{lemm:polarization}, so we omit it.
		
		To finish this section, we present a system of defining equations for the orbit $\Omega_{D^1,\xi}$. Namely, $\mu\in\ut^*$ belongs to $\Omega_{D^1,\xi}$ if and only if
		\begin{equation*}
			\begin{split}
				\Delta_{D^1}^{r,s}(\mu)&=\Delta_{D^1}^{r,s}(f_{D^1,\xi})\text{ for all }(r,s)\in(D^1\setminus R_1)\cup\{(1,n),(n-1,n)\},\\
				\gamma(\mu)&=\sum_{k=2}^{n-1}\mu_{1,k}\mu_{k,n}=\gamma(f_{D^1,\xi}).
			\end{split}
	\end{equation*}}
	
	The characters corresponding to 0-regular and 1-regular orbits (i.e., characters of depths 0 and~1) are computed at the end of the next section.

	\sect{The Mackey's little group method}\label{sect:Mackey_method}
	
	This section is devoted to the Mackey's method, which allows to reduce the study of representations of a semidirect product to the representations of its factors. This will allow us to reduce a description of characters of depth 2 to characters of depth 1.
	
	Let $G$ be a finite group represented as a semidirect product of its subgroups $A$ and $B$: $G=A\rtimes B$. This means that $A$ is a normal subgroup of $G$, $G=AB$ and $A\cap B=\{e\}$, where $e$ is a neutral element. It is evident that each element $g\in G$ can be uniquely expressed as a product $g=g_Ag_B$ for certain $g_A\in A$, $g_B\in B$. We will write $\pi^G_A(g)=g_A$, $\pi^G_B(g)=g_B$. Note that $\pi^G_B$ is a group homomorphism, while $\pi^G_A$ is not in general. Obviously, if $C$ is a subgroup of $B$ then $AC=A\rtimes C$ and $\pi^G_A(g)=\pi^{AC}_A(g)$, $\pi^G_B(g)=\pi^{AC}_C(g)$ for any element $g\in AC$.
	
	We will denote by $\Irr{K}$ the set of irreducible complex characters of a finite group $K$. Assume from now on that $A$ is abelian. Given a character $\kappa\in\Irr{A}$ (which is necessarily of degree 1, because $A$ is abelian) and an element $b\in B$, we define the character $\kappa^b\in\Irr{A}$ by the formula ${\kappa^b(a)=\kappa(bab^{-1})}$, $a\in A$. By definition, the \emph{centralizer} (or the \emph{little group}) of $\kappa$ in $B$ is the subgroup $$B^{\kappa}=\{b\in B\mid\kappa^b=\kappa\}.$$ Evidently, $\kappa_0=\kappa\circ\pi^{AB^{\kappa}}_A$ is an (irreducible) character of $AB^{\kappa}$ of degree 1.
	
	Next, pick a character $\psi\in\Irr{B^{\kappa}}$ and denote $\psi_0=\psi\circ\pi^{AB^{\kappa}}_{B^{\kappa}}\colon AB^{\kappa}\to\Cp$. Since $\pi^{AB^{\kappa}}_{B^{\kappa}}$ is a homomorphism, $\psi_0$ is again a character of $AB^{\kappa}$, and an elementary calculation shows that $\psi_0$ is in fact irreducible. Further, one can immediately check that $\kappa_0\psi_0\in\Irr{AB^{\kappa}}$. The Mackey's method can be formulated as follows (see, e.g., \cite[Proposition 1.3]{Lehrer74}).
	
	\mtheo{Let $G$\textup, $A$\textup, $B$ etc. be as above. Then \textup{i)} the induced character $\chi_{\kappa,\psi}=\Ind{AB^{\kappa}}{G}{\kappa_0\psi_0}$ of $G$ is irreducible\textup; \textup{ii)} each $\chi\in\Irr{G}$ has the form $\chi=\chi_{\kappa,\psi}$ for certain $\kappa\in\Irr{A}$\textup, $\psi\in\Irr{B^{\kappa}}$.}
	
	As an example of using of this method, in the rest of the section we will consider the case\break $G=U=U_n(\Fp_q)$. Let $\langle\cdot\rangle_{\Fp}$ denote the linear span over a field $\Fp$. Put
	\begin{equation*}
		\begin{split}
			A&=U'=\{g\in U\mid g_{i,j}=0\text{ for all }2\leq j<i\leq n\}=\exp\langle e_{i,1},~(i,1)\in\Co_1\rangle_{\Fp_q},\\
			B&=\wt U=\{g\in U\mid g_{i,1}=0\text{ for all }2\leq i\leq n\}=\exp\langle e_{i,j},~(i,j)\in\Phi(n)\setminus\Co_1\rangle_{\Fp_q}.
		\end{split}
	\end{equation*}
	It is clear that $U'\cong\Fp_q^{n-1}$ is abelian, $\wt U\cong U_{n-1}(\Fp_q)$ and $U=U'\rtimes\wt U$. Let $\chi\in\Irr{U}$ be a character of depth $k=0$ or 1, then it corresponds to the orbit $\Omega_{D^k,\xi}$ for a unique map $\xi\colon D^k\to\Fp_q^{\times}$. Our first goal is to present $\kappa\in\Irr{U'}$ and $\psi\in\Irr{\wt U^{\kappa}}$ such that $\chi=\chi_{\kappa,\psi}$. To do this, denote $f=f_{D^k,\xi}$ and set $\kappa=\theta\circ f\circ\ln$. In fact, for $a\in U'$,
	\begin{equation*}
		\kappa(a)=\theta(\xi(n-k,1)a_{n-k,1}).
	\end{equation*}
	
	One can easily check that the little group of $\kappa$ has the form
	\begin{equation*}
		\wt U^{\kappa}=\{b\in\wt U\mid b_{n-k,j}=0\text{ for all }2\leq j\leq n-k-1\}\cong
		\begin{cases}
			U_{n-2}(\Fp_q)&\text{ for }k=0,\\
			U_{n-2}(\Fp_q)\times\Fp_q&\text{ for }k=1.\\
		\end{cases}
	\end{equation*}
	Put also $V=\wt U^{\kappa}$ for $k=0$ and $V=\{b\in\wt U^{\kappa}\mid b_{n,n-1}=0\}$ for $k=1$, so that $V\cong U_{n-2}(\Fp_q)$. Let $\vt\cong\ut_{n-2}(\Fp_q)$ be the Lie algebra of~$V$, and set $\wt\Phi=\{(r,s)\in\Phi(n)\mid e_{r,s}\in\vt\}$. Denote by $\rho\colon\wt\Phi\to\Phi(n-2)$ the obvious bijection obtained by renumerating the indices of rows and columns, and put
	\begin{equation*}
		\wt D=\rho(D\cap\wt\Phi),~ \wt\xi=\restr{\xi}{D\cap\wt\Phi}\circ\rho^{-1}\colon\wt D\to\Fp_q^{\times},~ \wt f=f_{\wt D,\wt\xi}\in\ut_{n-2}(\Fp_q)^*.
	\end{equation*}
	Then the $U_{n-2}(\Fp_q)$-orbit of $\wt f$ is 0-regular.
	
	Note also that the set
	\begin{equation*}
		\Tu=\{g\in\wt U\mid g_{i,j}=0,\text{ if }i\neq n-k,~j<i\}=\exp\left(\langle e_{n-k,j},~(n-k,j)\in\Ro_{n-k}\setminus\Co_1\rangle_{\Fp_q}\right)
	\end{equation*}
	is a complete set or representatives for the left cosets $U/U'\wt U^{\kappa}$. Indeed, this follows directly from the well-known fact that each element of $U$ can be uniquely represented  as a product $$\prod_{(i,j)\in\Phi(n)}\exp(t_{i,j}e_{i,j}),$$ where $t_{i,j}\in\Fp_q$, and the product is taken in any fixed order. (In fact, $\exp(t_{i,j}e_{i,j})=1_n+t_{i,j}e_{i,j}$, where $1_n$ is the unit matrix.)
	
	\exam{Below\label{exam:0_1_decomposition_Mackey} we drew schematically the sets under consideration in the setup of Example~\ref{exam:0_1_regular_orbits}. The left and the right pictures correspond to the cases (i) $k=0$ and (ii) $k=1$ respectively. Boxes from $\wt\Phi$ are white (they correspond to $V$); boxes from $\Co_1$ are gray (they correspond to $U'$); boxes from $\Ro_{n-k}\setminus\Co_1$ are yellow (they correspond to $\Tu$); the box $(n,n-1)$ for $k=1$ is green (it corresponds to the factor $\Fp_q$ in the group $\wt U^{\kappa}\cong U_{n-2}(\Fp_q)\times\Fp_q$). Note that roots from $D\cap\wt\Phi=\rho^{-1}(\wt D)$ are marked by white symbols $\otimes$. We hope that it is clear now that the orbit of $\wt f$ is 0-regular in $\ut_{n-2}(\Fp_q)^*$.
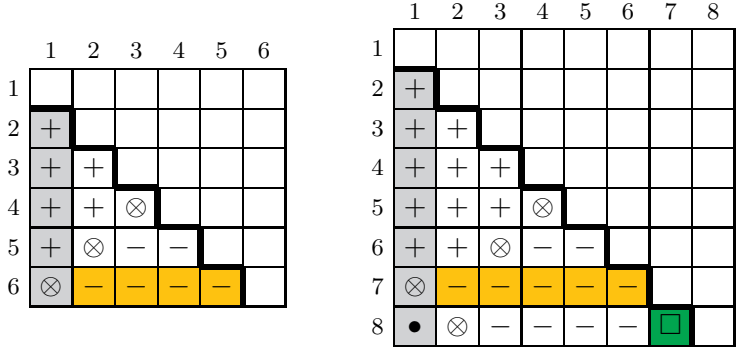
		\begin{equation*}
			\mymatrix{
				\Bot{2pt}\pho& \pho& \pho& \pho& \pho& \pho\\
				\gray\Rt{2pt}+& \Bot{2pt}\pho& \pho& \pho& \pho& \pho\\
				\gray+& \Rt{2pt}+& \Bot{2pt}\pho& \pho& \pho& \pho\\
				\gray+& +& \Rt{2pt}\otimes& \Bot{2pt}\pho& \pho& \pho\\
				\gray+& \otimes& -& \Rt{2pt}-& \Bot{2pt}\pho& \pho\\
				\gray\otimes& \yellow-& \yellow-& \yellow-& \yellow\Rt{2pt}-& \pho\\}\hspace{1cm}
			\mymatrix{
				\Bot{2pt}\pho& \pho& \pho& \pho& \pho& \pho& \pho& \pho\\
				\gray\Rt{2pt}+& \Bot{2pt}\pho& \pho& \pho& \pho& \pho& \pho& \pho\\
				\gray+& \Rt{2pt}+& \Bot{2pt}\pho& \pho& \pho& \pho& \pho& \pho\\
				\gray+& +& \Rt{2pt}+& \Bot{2pt}\pho& \pho& \pho& \pho& \pho\\
				\gray+& +& +& \Rt{2pt}\otimes& \Bot{2pt}\pho& \pho& \pho& \pho\\
				\gray+& +& \otimes& -& \Rt{2pt}-& \Bot{2pt}\pho& \pho& \pho\\
				\gray\otimes& \yellow-& \yellow-& \yellow-& \yellow-& \yellow\Rt{2pt}-& \Bot{2pt}\pho& \pho\\
				\gray\bullet& \otimes& -& -& -& -& \green\Rt{2pt}\square& \pho\\}
		\end{equation*}
	}
	
	For $k=0$, let $\psi$ be the irreducible character of $\wt U^{\kappa}=V\cong U_{n-2}(\Fp_q)$ corresponding to the 0-regular orbit $\Omega_{\wt D,\wt\xi}\subset\ut_{n-2}(\Fp_q)^*$. For $k=1$, we denote by $\psi$ the character of $\wt U^{\kappa}=V\times\Fp_q$ defined as follows. Let $\wt\psi$ be the irreducible character of $V\cong U_{n-2}(\Fp_q)$ corresponding to the 0-regular orbit $\Omega_{\wt D,\wt\xi}$. Notice that, given $g\in\wt U^{\kappa}$, the entries of the projection $g_V$ coincide (up to the obvious renumeration of indices) with the corresponding entries of $g$. We set $$\psi(g)=\wt\psi(g_V)\theta(\xi(n,n-1)g_{n,n-1}),$$ where $\xi(n,n-1)=0$ if $(n,n-1)\notin D^1$. The following proposition is one of the key steps in the obtaining an explicit formula for the character $\chi$.
	
	\mprop{One\label{prop:chi_0_1_regular_Mackey} has $\chi=\chi_{\kappa,\psi}$.}
	
	The proof of this proposition is similar to the proof of Proposition~\ref{prop:chi_2_regular_Mackey}, so we skip it. But the main ideas of the proof are as follows. First, recall the polarization $\pt$ for $f$ from Example~\ref{exam:0_1_regular_orbits}. One can easily check that $\pt\cap\vt$ is a polarization for $\wt f$ (after the identification $\vt\cong\ut_{n-2}(\Fp_q)$). Then, denote $\wt P=\exp(\pt\cap\vt)$ and recall the notion $P=\exp{\pt}$. It is easy to see that $P=\Tu\wt P$. It remains to use formula~(\ref{formula:chi_Omega}) and standard formulas for induced characters.
	
	To complete the calculation of $\chi$, we need to introduce some additional notation. Given a root $(i,j)\in\Phi(n)$, we set $\So_{i,j}=\{(i,l),~1\leq l<i\}\cup\{(l,j),~j<l\leq n\}$ (cf. the notion of the set $S(i,j)$ of singular roots). To a subset $D\subset\Phi(n)$ we attach the sets
	\begin{equation}
		\So_D=\bigcup_{(i,j)\in D}\So_{i,j},~\Ro_D=\Phi(n)\setminus\So_D.\label{formula:S_D_R_D}
	\end{equation}
	We will refer to the roots from $\Ro_D$ as $D$-\emph{regular} roots.
	
	\defi{A set $D\subset\Phi(n)$ is called 0-\emph{regular} if it is contained in $D^0$. It is called 1-\emph{regular} if $D'=D\setminus(\Co_1\cup\Ro_n)$ is a 0-regular subset not containing the root $(n-1,2)$, and $D''=D\cap(\Co_1\cap\Ro_n)$ is one of the following sets:
		\begin{equation*}
			\begin{split}
				&\varnothing,~\{(n,1)\},~\{(n-1,1)\},~\{(n,2)\},~\{(n-1,1),~(n,2)\},\\
				&\{(2,1),~(n,n-1)\},~\{(2,1),~(n-1,1),~(n,2),~(n,n-1)\}.
			\end{split}
		\end{equation*}
		(Note that in \cite{Ignatev09}, 0-regular and 1-regular subsets are called regular and 1-subregular respectively.)}
	
	\exam{Let\label{exam:1_subreg_subsets} $n=6$ or $n=8$. On the picture below we drew schematically 1-regular subsets\break $D=\{(5,1),~(6,2),~(4,3)\}\subset\Phi(6)$ and $D=\{(2,1),~(7,1),~(8,2),~(5,4),~(8,7)\}\subset\Phi(8)$. Roots from~$D$ are marked by $\otimes$, roots from $\Ro_D$ (i.e., $D$-regular roots) are blue. We also marked by red roots $(i,j)\in\Phi(n)$ for which there are no $(r,s)\in D$ such that $(i,j)\geq(r,s)$; their importance will be explained below.
		\begin{equation*}
			\mymatrix{
				\Bot{2pt}\pho& \pho& \pho& \pho& \pho& \pho\\
				\Rt{2pt}\red\pho& \Bot{2pt}\pho& \pho& \pho& \pho& \pho\\
				\red\pho& \red\Rt{2pt}\pho& \Bot{2pt}\pho& \pho& \pho& \pho\\
				\blue\pho& \blue\pho& \Rt{2pt}\otimes& \Bot{2pt}\pho& \pho& \pho\\
				\otimes& \pho& \blue\pho& \red\Rt{2pt}\pho& \Bot{2pt}\pho& \pho\\
				\blue\pho& \otimes& \blue\pho& \red\pho& \red\Rt{2pt}\pho& \pho\\}\hspace{1cm}
			\mymatrix{
				\Bot{2pt}\pho& \pho& \pho& \pho& \pho& \pho& \pho& \pho\\
				\Rt{2pt}\otimes& \Bot{2pt}\pho& \pho& \pho& \pho& \pho& \pho& \pho\\
				\blue\pho& \red\Rt{2pt}\pho& \Bot{2pt}\pho& \pho& \pho& \pho& \pho& \pho\\
				\blue\pho& \red\pho& \red\Rt{2pt}\pho& \Bot{2pt}\pho& \pho& \pho& \pho& \pho\\
				\blue\pho& \blue\pho& \blue\pho& \Rt{2pt}\otimes& \Bot{2pt}\pho& \pho& \pho& \pho\\
				\blue\pho& \pho& \pho& \blue\pho& \red\Rt{2pt}\pho& \Bot{2pt}\pho& \pho& \pho\\
				\blue\otimes& \pho& \pho& \blue\pho& \red\pho& \red\Rt{2pt}\pho& \Bot{2pt}\pho& \pho\\
				\blue\pho& \blue\otimes& \blue\pho& \blue\pho& \blue\pho& \blue\pho& \Rt{2pt}\otimes& \pho\\}
		\end{equation*}
	}
	We are ready now to describe conjugacy classes contained in the support of a character $\chi$ of depth $k=0$ or $1$. Namely, let $D\subset\Phi(n)$ be a $k$-regular subset, and $\vfi\colon D\to\Fp_q^{\times}$ be a map. Consider the element 
	\begin{equation*}
		x_D(\vfi)=1_n+\sum_{(i,j)\in D}\vfi(i,j)e_{i,j}\in U.\label{formula:x_D_vfi}
	\end{equation*}Let $\Ko_D(\vfi)\subset U$ be the conjugacy class of $x_D(\vfi)$.
	
	For a root $(r,s)\in\Ro_D$, let $$D^{(r,s)}=\{(i,j)\in D\mid(i,j)<(r,s)\}.$$ Denote $D^{(r,s)}\cup\{(r,s)\}=\{(i_1,j_1),\ldots,(i_t,j_t)\}$, where $j_1<\ldots<j_t$. For a matrix $x$, we denote by $M^D_{r,s}(x)$ the minor\label{deno:M_D_ij} of the matrix $x$ with the set of columns $j_1,\ldots,j_t$ and the set of rows $\sigma(i_1),\ldots,\sigma(i_t)$, where $\sigma$ is the unique permutation such that $\sigma(i_1)<\ldots<\sigma(i_t)$. (Cf. the notion of $D_{r,s}$ and $\Delta_D^{r,s}$ from Section~\ref{sect:orbit_method}.) Given $(i,j)\in\Phi(n)$, put also
	\begin{equation}
		\alpha_j(x)=\sum_{l=j+1}^{n-1}x_{n,l}x_{l,j},~\beta_i(x)=\sum_{l=2}^{i-1}x_{i,l}x_{l,1},~\gamma(x)=\sum_{l=2}^{n-1}x_{n,l}x_{l,1}=\alpha_1(x)=\beta_n(x).\label{formula:alpha_j_beta_i}
	\end{equation}
	
	\mprop{An\label{prop:0_1_depth_support} element $x\in U$ belongs to $\Ko_D(\vfi)$ if and only if $M_{i,j}^D(x)=M_{i,j}^D(x_D(\vfi))$ and\textup, for $k=1$ and $\{(2,1),~(n,n-1)\}\subset D$\textup, $\gamma(x)=\gamma(x_D(\vfi)),~\alpha_j(x)=\beta_i(x)=0\text{ for }(i,j)\in D\setminus(\Ro_n\cup\Co_1)$.
	}
	
	The case $k=0$ follows from \cite[Theorem 2.5, Proposition 3.1]{Andre01}, while the proof for the case $k=1$ can be found in \cite{Ignatev09}. It is essentially algebraic-geometric. Namely, it is proved that the ideal $J$ of $\Kp[U(\Kp)]$ generated by these polynomials is prime and that the corresponding affine variety $V(J)$ of its common zeroes is $U(\Kp)$-invariant. Then it is checked that the centralizer of $x_D(\vfi)$ in $U(\Kp)$ is an affine subvariety isomorphic to the affine space of dimension $\codim_{U(\Kp)}V(J)$. Thus, it remains to note that $x_D(\vfi)\in V(J)$ and restrict to $\Fp_q$-points. We use a similar scheme to describe conjugacy classes in Section~\ref{sect:main_result} below.
	
	\exam{The conjugacy class $\Ko_D(\vfi)$ corresponding to the right picture from Example \ref{exam:1_subreg_subsets} is defined by the equations $x_{i,j}=0$ if the root $(i,j)$ is red, $x_{4,3}=\vfi(4,3)$, and
		\begin{equation*}
			\begin{vmatrix}
				x_{4,j}&x_{4,3}\\
				x_{i,j}&x_{i,3}\\
			\end{vmatrix}
			=-\vfi(i,j)\vfi(4,3)\text{ for }(i,j)\in\{(5,1),~(6,2)\}.
		\end{equation*}
		At the same time, $\Ko_D(\vfi)$ corresponding to the right picture is defined by the equations $x_{i,j}=0$ if the root $(i,j)$ is red, $x_{i,j}=\vfi(i,j)$ for $(i,j)\in\{(2,1),~(5,4),~(8,7)\}$, and
		\begin{equation*}
			\begin{split}
				&\begin{vmatrix}
					x_{5,j}&x_{5,4}\\
					x_{i,j}&x_{i,4}\\
				\end{vmatrix}=0\text{ for }i\in\{6,7\},~j\in\{2,3\},~\gamma(x)=\sum_{l=2}^7x_{8,l}x_{l,1}=\vfi(8,2)\vfi(2,1)+\vfi(8,7)\vfi(7,1),\\
				&\alpha_4(x)=x_{8,5}x_{5,4}+x_{8,6}x_{6,4}+x_{8,7}x_{7,4}=\beta_5(x)=x_{5,2}x_{2,1}+x_{5,3}x_{3,1}+x_{5,4}x_{4,1}=0.
			\end{split}
		\end{equation*}
		One can see that each white or red root gives us an equation, plus some extra equations for the case $k=1$, $\{(2,1),~(n,n-1)\}\subset D$.}
	
	Finally, to formulate an answer, we will define certain integer numbers. Recall from Section~\ref{sect:orbit_method} that $\Mo=\{(i,j)\in\Phi(n)\mid i>n-j+1\}$. Given a $k$-regular subset $D\subset\Phi(n)$ for $k=0$ or 1, we set
	\begin{equation}
		m_D=\begin{cases}
			|\Ro_D\cap\Mo|-k,&\text{if }\{(2,1),(n,n-1)\}\not\subset D,\\
			|\Ro_D\cap\Mo|+n-1,&\text{if }\{(2,1),(n,n-1)\}\subset D.\\
		\end{cases}\label{formula:m_D_0_1}
	\end{equation} The following theorem gives us an explicit description of the support of $\chi$, as well as an explicit formula how to compute the value $\chi$ on an element from its support. Given a $k$-regular orbit $\Omega_{D^k,\xi}$, recall the notion of the linear form $f=f_{D^k,\xi}$, $k=0$ or 1.
	
	\mtheo{Let\label{theo:0_1_depth_character} $\chi$ be the character of depth $k$ corresponding to a $k$-regular orbit $\Omega_{D^k,\xi}$. Then \textup{i)} $\Supp{\chi}=\bigcup\Ko_{D,\xi}$, where the union is taken over all $k$-regular subsets $D\subset\Phi(n)$ and all maps $\xi\colon D\to\Fp_q^{\times}$ satisfying the additional condition $\xi(n-1,1)\vfi(2,1)=\xi(n,2)\vfi(n,n-1)$ if $k=1$.\break \textup{ii)} If $g$ belongs to $\Ko_D(\vfi)$ satisfying \textup{(i)} then
		\begin{equation*}
			\chi(g)=\chi(x_D(\vfi))=q^{m_D}\theta(f(x_D(\vfi)-1_n))=q^{m_D}\prod_{(i,j)\in D\cap D^k}\theta(\xi(i,j)\vfi(i,j)).
	\end{equation*}}
	
	The proof can be found in \cite[Section 5]{Ignatev09}. It uses induction by $n$ and Proposition~\ref{prop:chi_0_1_regular_Mackey}. The proof of our main result presented in Section~\ref{sect:proofs_classes} follows the same strategy, but the calculations are much more technical.
	
	\sect{Classification of 2-regular orbits}\label{sect:2_orbits}
	
	From now on and to the end of the paper we assume that $n\geq5$. Before calculating characters of depth 2 we need to classify all 2-regular orbits, of course. To do this, recall the notion of $D^0$ from Example~\ref{exam:0_1_regular_orbits} (i) and denote by $D^2$ the following subset:
	\begin{equation*}
		D^2=(D^0\setminus\{(n,1),~(n-1,2),~(n-2,3)\})\cup\{(n-1,1)~(n-2,2),~(n,3)\}\cup R_2,
	\end{equation*}
	where $R_2$ is a (possibly, empty) subset of $\{(n,n-2),~(n,n-1)\}$. As above, to each map $\xi\colon D^2\to\Fp_q^{\times}$, one can assign the linear form $$f_{D^2,\xi}=\sum\nolimits_{(i,j)\in D^2}\xi(i,j)e_{j,i}\in\ut^*,$$ and denote by $\Omega_{D^2,\xi}$ the coadjoint orbit of this linear form. As we will see below, the orbits of the form~$\Omega_{D^2,\xi}$ are exactly the 2-regular orbits, i.e., the orbits of maximal possible dimension in the invariant subspace $\Xu_2$ of $\ut^*$.
	
	\exam{On\label{exam:2_regular_orbits} the picture below we drew the set $D^2$ for $n=10$. As in Example~\ref{exam:0_1_regular_orbits}, roots from $D^2$ are marked by $\otimes$, roots from $$\Mo=\left(\bigcup_{(i,j)\in D}S^-(i,j)\right)\setminus\{(n,n-2),~(n,n-1)\}$$ are marked by minuses, the roots $(n-1)$ and $(n,2)$ are marked by dots, the roots $(n,n-2)$ and $(n,n-1)$ are marked by $\square$'s, and all other roots, i.e., roots from $$\left(\bigcup_{(i,j)\in D}S^+(i,j)\right)\setminus\{(n-1,3),~(n-2,3)\},$$ are marked by pluses.
		\begin{equation*}
			\mymatrix{
				\Bot{2pt}\pho& \pho& \pho& \pho& \pho& \pho& \pho& \pho& \pho& \pho\\
				\Rt{2pt}+& \Bot{2pt}\pho& \pho& \pho& \pho& \pho& \pho& \pho& \pho& \pho\\
				+& \Rt{2pt}+& \Bot{2pt}\pho& \pho& \pho& \pho& \pho& \pho& \pho& \pho\\
				+& +& \Rt{2pt}+& \Bot{2pt}\pho& \pho& \pho& \pho& \pho& \pho& \pho\\
				+& +& +& \Rt{2pt}+& \Bot{2pt}\pho& \pho& \pho& \pho& \pho& \pho\\
				+& +& +& +& \Rt{2pt}\otimes& \Bot{2pt}\pho& \pho& \pho& \pho& \pho\\
				+& +& +& \otimes& -& \Rt{2pt}-& \Bot{2pt}\pho& \pho& \pho& \pho\\
				+& \otimes& -& -& -& -& \Rt{2pt}-& \Bot{2pt}\pho& \pho& \pho\\
				\otimes& -& -& -& -& -& -& \Rt{2pt}-& \Bot{2pt}\pho& \pho\\
				\bullet& \bullet& \otimes& -& -& -& -& \square& \Rt{2pt}\square& \pho\\}
	\end{equation*}}
	
	\lemmp{The\label{lemm:polarization} subspace $\pt=\langle e_{i,j},~(i,j)\in\Phi(n)\setminus\Mo\rangle_{\Fp_q}\subset\ut$ is a polarization for $f=f_{D^2,\xi}$.}{It is evident that $\pt$ is a subalgebra of $\ut$. Note that, given $i>l>j$, one has $[e_{i,l},e_{l,j}]=e_{i,j}$ if and only if $(i,l)\in S^-(i,j)$ and $(l,j)\in S^+(i,j)$. Note also that if $(i,j)\in D^2$ then $\Mo$ contains exactly one of the roots $(i,l),~(l,j)\in S(i,j)$ for each $l$ between $i$ and $j$. Since $\Supp{f}=D^2$, this implies that $\pt$ is $\beta_f$-isotropic subspace.
		
		It is easy to see that if a root $(i,l)\in\Phi(n)$ is marked by plus then there exists a unique root $(l,j)\in\Mo$ such that $(i,j)\in D^2$. Now, assume that $x=y+z\in\ut$ belongs to a $\beta_f$-isotropic space containing $\pt$ for $y\in\pt$, $z\in\langle e_{i.l},~(i,l)\in\Mo\rangle_{\Fp_q}$ and $z_{i,l}\neq0$ for some $(i,l)\in\Mo$. Then there exists $(l,j)\in\Phi(n)\setminus\Mo$ marked by plus such that $(i,j)\in D^2$. Hence, $$f([x,e_{l,j}])=\xi(i,j)z_{i,l}\neq0,$$ a contradiction. Thus, $\pt$ is a maximal $\beta_f$-isotropic subspace, as required.}
	
	\corop{An\label{coro:2_regular_orbit} orbit of the form $\Omega_f$ is $2$-regular.}{Clearly, $f\in\Xu_2$. But the maximal possible dimension of a $U(\Kp)$-orbit of a linear form from~$\Xu_2$ equals $N-4$, see the Introduction. At the same time, $\dim\Omega_f(\Kp)=2\codim_{\ut}\pt=N-4$. This completes the proof.}
	
	\nota{Using\label{nota:basis_defi_ideal} the methods of \cite[Theorem 1.9]{IgnatevPanov09}, one can construct a system of defining equations for the orbit $\Omega_f$ (of $\Omega_f(\Kp)$). Namely, let $\succ$ be the following linear order on $\Phi(n)$ refining the usual order $\geq$: $(a,b)\succ(c,d)$ if $b<d$ or $b=d$, $a>c$. To each $(i,j)$ marked by $\otimes$, $\square$ or $\bullet$ one can attach a polynomial $Q_{i,j}$ in the symmetric algebra $S(\ut)$ such that i) the defining ideal of $\Omega_f(\Kp)$ is generated by $Q_{i,j}-Q_{i,j}(f)$ for all such $(i,j)$; ii) each polynomial has the form $Q_{i,j}=\lambda_{i,j}+P_{i,j}^{\succ}$. Here $\lambda_{r,s}$ are the standard coordinate functions and $P_{i,j}^{\succ}$ belongs to the subalgebra of $S(\ut)$ generated by $\lambda_{r,s}$ for $(r,s)\succ(i,j)$ \cite[Corollary 1.11]{IgnatevPanov09}. Of course, this immediately implies again that $\dim\Omega_f(\Kp)$ equals the number of pluses and minuses on the picture, i.e., $N-4$. Furthermore, for $(i,j)$ marked by $\otimes$, one has $Q_{i,j}=\Delta_{D^2}^{i,j}$.}
	
	We introduce some more notation. Let $A$, $B$, $C$ be the following subvarieties of $\Xu_2(\Kp)=\Xu_2\otimes_{\Fp_q}\Kp$:
	\begin{equation*}
		\begin{split}
			A&=\{\lambda\in\Xu_2(\Kp)\mid\lambda_{n,3}=0\},\\
			B&=\{\lambda\in\Xu_2(\Kp)\mid\lambda_{n-1,1}=0\},\\
			C&=\left\{\lambda\in\Xu_2(\Kp)\mid\begin{vmatrix}\lambda_{n-2,1}&\lambda_{n-2,2}\\
				\lambda_{n-1,1}&\lambda_{n-1,2}\end{vmatrix}=0\right\}.
		\end{split}
	\end{equation*}
	Recall that one can consider $\ut^*$ as an $\Fp_q$-subspace of $\ut^*(\Kp)$.
	
	\lemmp{Let\label{lemm:not_2_regular} $\lambda\in\Xu_2$. Assume that $\lambda\in A\cup B\cup C$. Then the orbit of $\lambda$ is not $2$-regular.}{It is easy to check that the subvarieties $A$, $B$, $C$ are $U(\Kp)$-invariant. Consider the following subsets of $\Phi(n)$:
		\begin{equation*}
			\begin{split}
				D_A&=(D^2\setminus\{(n,3),~(n-3,4)\})\cup\{(n-3,3),~(n,4),~(n,n-3)\},\\
				D_B&=(D^2\setminus\{(n-1,1),~(n-2,2)\})\cup\{(n-2,1),~(n-1,2),~(n-1,n-2)\},\\
				D_C&=(D^2\setminus\{(n-2,2),~(n-3,4),~(n,n-2)\})\cup\{(n-3,2),~(n-2,4),~(n-2,n-3)\}.
			\end{split}
		\end{equation*}
		On the picture below, we drew schematically these subsets for $n=10$ ($D_A$ is on the left, $D_B$ is in the middle, and $D_C$ is on the right). We use marks similar to the previous example.
		\begin{equation*}
			\hspace{-0.8cm}\mymatrix{
				\Bot{2pt}\pho& \pho& \pho& \pho& \pho& \pho& \pho& \pho& \pho& \pho\\
				\Rt{2pt}+& \Bot{2pt}\pho& \pho& \pho& \pho& \pho& \pho& \pho& \pho& \pho\\
				+& \Rt{2pt}+& \Bot{2pt}\pho& \pho& \pho& \pho& \pho& \pho& \pho& \pho\\
				+& +& \Rt{2pt}+& \Bot{2pt}\pho& \pho& \pho& \pho& \pho& \pho& \pho\\
				+& +& +& \Rt{2pt}+& \Bot{2pt}\pho& \pho& \pho& \pho& \pho& \pho\\
				+& +& +& +& \Rt{2pt}\otimes& \Bot{2pt}\pho& \pho& \pho& \pho& \pho\\
				+& +& \otimes& -& -& \Rt{2pt}-& \Bot{2pt}\pho& \pho& \pho& \pho\\
				+& \otimes& -& -& -& -& \Rt{2pt}-& \Bot{2pt}\pho& \pho& \pho\\
				\otimes& -& -& -& -& -& -& \Rt{2pt}-& \Bot{2pt}\pho& \pho\\
				\bullet& \bullet& \bullet& \otimes& -& -& \square& \square& \Rt{2pt}\square& \pho\\}~   
			\mymatrix{
				\Bot{2pt}\pho& \pho& \pho& \pho& \pho& \pho& \pho& \pho& \pho& \pho\\
				\Rt{2pt}+& \Bot{2pt}\pho& \pho& \pho& \pho& \pho& \pho& \pho& \pho& \pho\\
				+& \Rt{2pt}+& \Bot{2pt}\pho& \pho& \pho& \pho& \pho& \pho& \pho& \pho\\
				+& +& \Rt{2pt}+& \Bot{2pt}\pho& \pho& \pho& \pho& \pho& \pho& \pho\\
				+& +& +& \Rt{2pt}+& \Bot{2pt}\pho& \pho& \pho& \pho& \pho& \pho\\
				+& +& +& +& \Rt{2pt}\otimes& \Bot{2pt}\pho& \pho& \pho& \pho& \pho\\
				+& +& +& \otimes& -& \Rt{2pt}-& \Bot{2pt}\pho& \pho& \pho& \pho\\
				\otimes& -& -& -& -& -& \Rt{2pt}-& \Bot{2pt}\pho& \pho& \pho\\
				\bullet& \otimes& -& -& -& -& -& \Rt{2pt}\square& \Bot{2pt}\pho& \pho\\
				\bullet& \bullet& \otimes& -& -& -& -& \square& \Rt{2pt}\square& \pho\\}~
			\mymatrix{
				\Bot{2pt}\pho& \pho& \pho& \pho& \pho& \pho& \pho& \pho& \pho& \pho\\
				\Rt{2pt}+& \Bot{2pt}\pho& \pho& \pho& \pho& \pho& \pho& \pho& \pho& \pho\\
				+& \Rt{2pt}+& \Bot{2pt}\pho& \pho& \pho& \pho& \pho& \pho& \pho& \pho\\
				+& +& \Rt{2pt}+& \Bot{2pt}\pho& \pho& \pho& \pho& \pho& \pho& \pho\\
				+& +& +& \Rt{2pt}+& \Bot{2pt}\pho& \pho& \pho& \pho& \pho& \pho\\
				+& +& +& +& \Rt{2pt}\otimes& \Bot{2pt}\pho& \pho& \pho& \pho& \pho\\
				+& \otimes& -& -& -& \Rt{2pt}-& \Bot{2pt}\pho& \pho& \pho& \pho\\
				+& \bullet& +& \otimes& -& -& \Rt{2pt}\square& \Bot{2pt}\pho& \pho& \pho\\
				\otimes& -& -& -& -& -& -& \Rt{2pt}-& \Bot{2pt}\pho& \pho\\
				\bullet& \bullet& \otimes& -& -& -& \square& -& \Rt{2pt}\square& \pho\\}
		\end{equation*}
		
		Let $D$ be one of the subsets $D_A$, $D_B$, $D_C$, and $\xi\colon D\to\Fp_q$ be a map such that $\xi(i,j)\neq0$ if $(i,j)$ is marked by $\otimes$; put $\lambda=f_{D,\xi}$. We claim that $\dim\Omega_{\lambda}(\Kp)$ coincides with the total number of pluses and minuses. (Note that the number of pluses equals the number of minuses for each picture, and the total number of pluses and minuses equals $N-6$.) Indeed, according to \cite[Corollary 1.11]{IgnatevPanov09}, the defining ideal of this orbit has a basis completely similar to the one described in Remark~\ref{nota:basis_defi_ideal}. Equivalently, let $\Mo$ be the set of roots marked by minuses. One can easily prove in the spirit of Lemma~\ref{lemm:polarization} that the subspace $$\pt=\langle e_{i,j},~(i,j)\in\Phi(n)\setminus\Mo\rangle_{\Fp_q}\subset\ut$$ is a polarization for $\lambda$.
		
		Now, consider $\bigcup\Omega_{\lambda}$, where the union is taken over all maps $\xi\colon D\to\Kp$ satisfying $\xi(i,j)\neq0$ for $(i,j)$ marked by $\otimes$. As it was mentioned in Remark~\ref{nota:basis_defi_ideal}, this union has the form $$\{\mu\in\ut^*(\Kp)\mid\Delta_D^{i,j}(\mu)\neq0,~(i,j)\text{ is marked by }\otimes\}.$$ Hence, this union is open and dense in the corresponding variety $A$, $B$ or $C$. But the set of linear forms having the orbit of maximal dimension is also dense~\cite{Kraft00}, thus the maximal possible dimension of an orbit in this subvariety is $N-6$. The proof is complete.}
	
	Note that the fact that the maximal possible dimension of an orbit in $A$ or $B$ is $N-6$ follows directly from Panov's results \cite[Theorem 2]{Panov09}.
	
	\theop{\textup{i)}The orbit of each\label{theo:2_regular_orbits} linear form $f_{D^2,\xi}$ is \textup2-regular. \textup{ii)} Each \textup2-regular orbit contains exactly one point of the form $f_{D^2,\xi}$ for a certain map $\xi\colon D^2\to\Fp_q^{\times}$.}{i) See Corollary~\ref{coro:2_regular_orbit}.
		
		ii) Assume that $\lambda\in\Xu_2$ has a 2-regular orbit. Lemma~\ref{lemm:not_2_regular} says that $\lambda\notin A\cup B\cup C$, i.e.,
		\begin{equation*}
			\lambda_{n-1,1}\neq0,~\lambda_{n,3}\neq0,~\begin{vmatrix}\lambda_{n-2,1}&\lambda_{n-2,2}\\\lambda_{n-1,1}&\lambda_{n-1,2}\end{vmatrix}\neq0.
		\end{equation*}
		
		In \cite{Andre95}, Andre established the following decomposition of $\ut^*(\Kp)$ into $U(\Kp)$-invariant subvarieties. Given a rook placement $D\subset\Phi(n)$ and a map $\xi\colon D\to\Kp^{\times}$, one can define the following $U(\Kp)$-invariant subvariety of $\ut^*(\Kp)$: $$\Ou_{D,\xi}(\Kp)=\{\mu\in\ut^*(\Kp)\mid\Delta_D^{i,j}(\mu)=\Delta_D^{r,s}(f_{D,\xi})\text{ for all }(r,s)\in\bigcup_{(i,j)\in D}S(i,i)\}.$$
		It turns out that $$\ut^*(\Kp)=\bigsqcup_{D,\xi}\Ou_{D,\xi}(\Kp).$$ Given a map $\xi\colon D\to\Fp_q^{\times}$, one can similarly define the $U$-invariant subset $\Ou_{D,\xi}\subset\ut^*$ and obtain the decomposition $\ut^*=\bigsqcup_{D,\xi}\Ou_{D,\xi}$. Clearly, $\Ou_{D,\xi}$ is the set of $\Fp_q$-points of $\Ou_{D,\xi}(\Kp)$.
		
		Hence, there exist the unique rook placement $D\subset\Phi(n)$ and the unique map $\xi\colon D\to\Fp_q^{\times}$ such that $\lambda\in\Ou_{D,\xi}$. Furthermore, it follows from the proof of \cite[Proposition 3]{Andre95} that $$\{(n-1,1),~(n-2,2),~(n,3)\}\subset D.$$ It is easy to check that if $D$ does not coincide with the set of all roots from $D^2$ marked by $\otimes$ (i.e., with $D^2\setminus R_2$) then $\dim\Ou_{D,\xi}(\Kp)<N-2$. Actually, this follows from the fact that $$\dim\Ou_{D,\xi}(\Kp)=\sum_{(i,j)\in D}|S(i,j)|,$$ see \cite[Theorem 2]{Andre95}.
		
		As it was mentioned in Remark~\ref{nota:basis_defi_ideal}, there exist $U(\Kp)$-invariant polynomials $Q_{n,n-1}$ and $Q_{n,n-2}$, which are clearly algebraically independent with all the minors $\Delta_D^{r,s}$ for $(r,s)\in\bigcup_{(i,j)\in D}S(i,j)$. Hence, given scalars $c_1$, $c_2\in\Kp$, the set
		$$Z_{c_1,c_2}(\Kp)=\{\mu\in\Ou_{D,\xi}(\Kp)\mid Q_{n,n-i}(\mu)=c_i,~i=1,2\}$$
		is a $U(\Kp)$-invariant subvariety of $\Ou_{D,\xi}(\Kp)$ of dimension $\dim\Ou_{D,\xi}-2$, and $$\Ou_{D,\xi}(\Kp)=\bigsqcup_{c_1,c_2}Z_{c_1,c_2}(\Kp).$$ Similarly, we have the decomposition $$\Ou_{D,\xi}=\bigsqcup\limits_{c_1,c_2} Z_{c_1,c_2},$$ where $Z_{c_1,c_2}$ is the set of $\Fp_q$-points of $Z_{c_1,c_2}(\Kp)$, and the union is taken over all $c_1$, $c_2\in\Fp_q$.
		
		Thus, the orbit $\Omega_{\lambda}(\Kp)$ is contained in $Z_{c_1,c_2}(\Kp)$ for unique $c_1$, $c_2\in\Fp_q$. Since $\dim\Omega_{\lambda}(\Kp)=N-2$, we conclude that $D=D^2\setminus R_2$. It is evident that there exist the unique subset $D^2$ and the unique map from $D^2$ to $\Fp_q^{\times}$ extending the map $\xi\colon D\to\Fp_q^{\times}$ (which we denote again by $\xi$) such that $c_i=Q_{n,n-i}(f_{D,\xi})$ for $i=1$, $2$. According to Remark~\ref{nota:basis_defi_ideal}, $Z_{c_1,c_2}$ has the same system of defining equation as $\Omega_{\lambda}$, so they coincide. This completes the proof.
	}
	
	Theorem~\ref{theo:2_regular_orbits} gives us a complete classification of 2-regular orbits on $\ut^*$. Thus, we can move to the our main result, an explicit description of the corresponding irreducible characters, i.e., the character of depth 2. The next section contains statements of the results, while the proofs are given in Sections~\ref{sect:proofs_classes} and \ref{sect:support} below.
	
	\sect{The main result}\label{sect:main_result}
	
	Let $\Omega_{D^2,\xi}=\Omega_{f_{D^2,\xi}}$ be a 2-regular orbit and $\chi=\chi_{D^2,\xi}$ be the corresponding character of depth~2. Our first goal in this section is to describe its support, i.e., to present a system of equations defining conjugacy classes which belong to $\Supp{\chi}$. To do this, we will introduce certain subsets of $\Phi(n)$ and subvarieties of $U$ as follows. Put
	\begin{equation*}
		\begin{split}
			D_0=&\{(n,1),(n,2),(n,3),(n-1,1),(n-2,2)\},\\
			D_2=&\{(3,2),(n,n-2)\}\text{, }\{(3,2),(n,n-2),(n-2,1)\}\text{, }\\
			&\{(3,2),(n,n-2),(n,3)\}\text{ or }\{(3,2),(n,n-2),(n-2,1),(n,3)\},\\
			D_3=&\{(3,1),(n,n-1)\}\text{, }\{(3,1),(n,n-1),(n-2,2)\},\\
			&\{(3,1),(n,n-1),(n,3)\}\text{ or }\{(3,1),(n,n-1),(n-2,2),(n,3)\},\\
			D_4=&\{(3,2),(n,n-1)\}\text{, }\{(3,2),(n,n-1),(n,3)\}\text{, }\\
			&\{(3,2),(n,n-1),(n-1,1)\}\text{ or }\{(3,2),(n,n-1),(n,3),(n-1,1)\},\\
			D_5=&\{(3,2),(n,n-1),(n-2,1),(n-1,2)\}\text{ or }\\
			&\{(3,2),(n,n-1),(n,3),(n-1,2),(n-2,1)\}.\\    
		\end{split}
	\end{equation*}
	
	\defi{We say that a subset $D\subset\Phi(n)$ is 2-\emph{regular} if $D=D'\cup D_i$, where $D'$ is a regular subset which does not contain roots $(n,1),(n-1,2),(n-2,3)$, $i=1,2,3,4,5$. Here $D_1$ is an arbitrary rook placement contained in $D_0$. Next, we will split an arbitrary subset $D\in\Phi(n)$ into the union
		$D=\overline{D}\sqcup D^*$, where $D^*=\{(i,j)\in D\mid i=3,n-1 \text{ or } n\}$ for convenience.}
	
	\exam{Below we schematically drew some examples of the set $D$ containing $D_1$, $D_2$, $D_3$, $D_4$ or $D_5$. Roots from $D$ are marked by $\otimes$, and roots from $D_i$ are colored in red. 
		\begin{equation*}
			D\supset D_1\colon\mymatrix{
				\Bot{2pt}\pho& \pho& \pho& \pho& \pho& \pho& \pho& \pho& \pho& \pho\\
				\Rt{2pt}& \Bot{2pt}\pho& \pho& \pho& \pho& \pho& \pho& \pho& \pho& \pho\\
				\pho& \Rt{2pt}& \Bot{2pt}\pho& \pho& \pho& \pho& \pho& \pho& \pho& \pho\\
				\pho& \pho& \Rt{2pt}& \Bot{2pt}\pho& \pho& \pho& \pho& \pho& \pho& \pho\\
				\pho& \pho& \pho& \Rt{2pt}& \Bot{2pt}\pho& \pho& \pho& \pho& \pho& \pho\\
				\pho& \pho& \pho& \pho& \Rt{2pt}\otimes& \Bot{2pt}\pho& \pho& \pho& \pho& \pho\\
				\pho& \pho& \pho& \otimes& \pho& \Rt{2pt}\pho& \Bot{2pt}\pho& \pho& \pho& \pho\\
				\pho& \red\otimes& \pho& \pho& \pho& \pho& \Rt{2pt}\pho& \Bot{2pt}\pho& \pho& \pho\\
				\red\otimes& \pho& \pho& \pho& \pho& \pho& \pho& \Rt{2pt}& \Bot{2pt}\pho& \pho\\
				\pho& \pho& \red\otimes& \pho& \pho& \pho& \pho& \pho& \Rt{2pt}& \pho\\}
			\hspace{1.2cm}D\supset D_2\colon\mymatrix{
				\Bot{2pt}\pho& \pho& \pho& \pho& \pho& \pho& \pho& \pho& \pho& \pho\\
				\Rt{2pt}& \Bot{2pt}\pho& \pho& \pho& \pho& \pho& \pho& \pho& \pho& \pho\\
				\pho& \red\Rt{2pt}\otimes& \Bot{2pt}\pho& \pho& \pho& \pho& \pho& \pho& \pho& \pho\\
				\pho& \pho& \Rt{2pt}& \Bot{2pt}\pho& \pho& \pho& \pho& \pho& \pho& \pho\\
				\pho& \pho& \pho& \Rt{2pt}& \Bot{2pt}\pho& \pho& \pho& \pho& \pho& \pho\\
				\pho& \pho& \pho& \pho& \Rt{2pt}\otimes& \Bot{2pt}\pho& \pho& \pho& \pho& \pho\\
				\pho& \pho& \pho& \otimes& \pho& \Rt{2pt}\pho& \Bot{2pt}\pho& \pho& \pho& \pho\\
				\pho& \pho& \pho& \pho& \pho& \pho& \Rt{2pt}\pho& \Bot{2pt}\pho& \pho& \pho\\
				\red\otimes& \pho& \pho& \pho& \pho& \pho& \pho& \Rt{2pt}& \Bot{2pt}\pho& \pho\\
				\pho& \pho& \red\otimes& \pho& \pho& \pho& \pho& \red\otimes& \Rt{2pt}& \pho\\}
		\end{equation*}
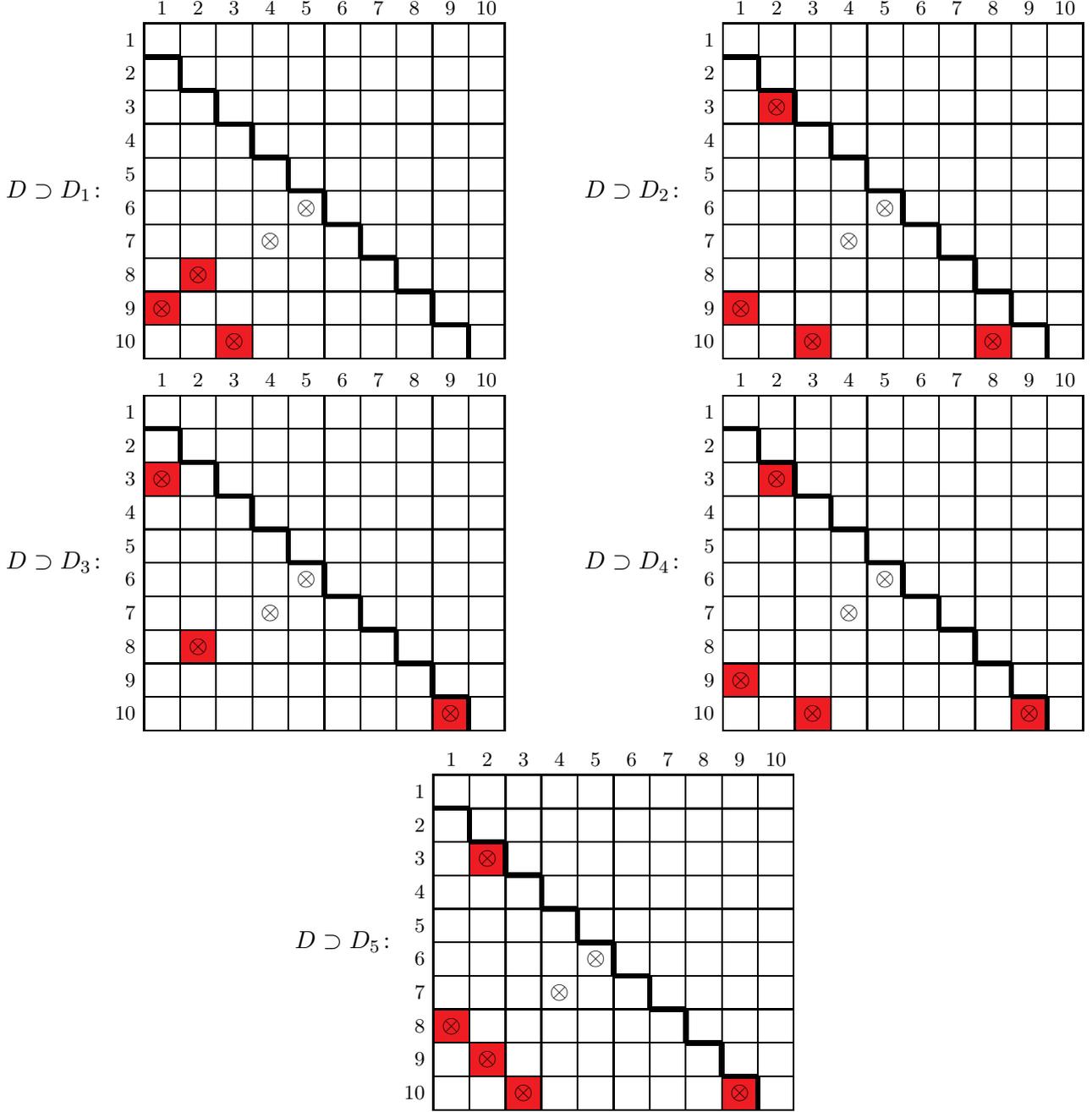
		\begin{equation*}
			D\supset D_3\colon\mymatrix{
				\Bot{2pt}\pho& \pho& \pho& \pho& \pho& \pho& \pho& \pho& \pho& \pho\\
				\Rt{2pt}& \Bot{2pt}\pho& \pho& \pho& \pho& \pho& \pho& \pho& \pho& \pho\\
				\red\otimes& \Rt{2pt}& \Bot{2pt}\pho& \pho& \pho& \pho& \pho& \pho& \pho& \pho\\
				\pho& \pho& \Rt{2pt}& \Bot{2pt}\pho& \pho& \pho& \pho& \pho& \pho& \pho\\
				\pho& \pho& \pho& \Rt{2pt}& \Bot{2pt}\pho& \pho& \pho& \pho& \pho& \pho\\
				\pho& \pho& \pho& \pho& \Rt{2pt}\otimes& \Bot{2pt}\pho& \pho& \pho& \pho& \pho\\
				\pho& \pho& \pho& \otimes& \pho& \Rt{2pt}\pho& \Bot{2pt}\pho& \pho& \pho& \pho\\
				\pho& \red\otimes& \pho& \pho& \pho& \pho& \Rt{2pt}\pho& \Bot{2pt}\pho& \pho& \pho\\
				\pho& \pho& \pho& \pho& \pho& \pho& \pho& \Rt{2pt}& \Bot{2pt}\pho& \pho\\
				\pho& \pho& \pho& \pho& \pho& \pho& \pho& \pho& \red\Rt{2pt}\otimes& \pho\\}
			\hspace{1.2cm}D\supset D_4\colon\mymatrix{
				\Bot{2pt}\pho& \pho& \pho& \pho& \pho& \pho& \pho& \pho& \pho& \pho\\
				\Rt{2pt}& \Bot{2pt}\pho& \pho& \pho& \pho& \pho& \pho& \pho& \pho& \pho\\
				\pho& \red\Rt{2pt}\otimes& \Bot{2pt}\pho& \pho& \pho& \pho& \pho& \pho& \pho& \pho\\
				\pho& \pho& \Rt{2pt}& \Bot{2pt}\pho& \pho& \pho& \pho& \pho& \pho& \pho\\
				\pho& \pho& \pho& \Rt{2pt}& \Bot{2pt}\pho& \pho& \pho& \pho& \pho& \pho\\
				\pho& \pho& \pho& \pho& \Rt{2pt}\otimes& \Bot{2pt}\pho& \pho& \pho& \pho& \pho\\
				\pho& \pho& \pho& \otimes& \pho& \Rt{2pt}\pho& \Bot{2pt}\pho& \pho& \pho& \pho\\
				\pho& \pho& \pho& \pho& \pho& \pho& \Rt{2pt}\pho& \Bot{2pt}\pho& \pho& \pho\\
				\red\otimes& \pho& \pho& \pho& \pho& \pho& \pho& \Rt{2pt}& \Bot{2pt}\pho& \pho\\
				\pho& \pho& \red\otimes& \pho& \pho& \pho& \pho& \pho& \red\Rt{2pt}\otimes& \pho\\}
		\end{equation*}
	}
	\begin{equation*}
		D\supset D_5\colon\mymatrix{
			\Bot{2pt}\pho& \pho& \pho& \pho& \pho& \pho& \pho& \pho& \pho& \pho\\
			\Rt{2pt}& \Bot{2pt}\pho& \pho& \pho& \pho& \pho& \pho& \pho& \pho& \pho\\
			\pho& \red\Rt{2pt}\otimes& \Bot{2pt}\pho& \pho& \pho& \pho& \pho& \pho& \pho& \pho\\
			\pho& \pho& \Rt{2pt}& \Bot{2pt}\pho& \pho& \pho& \pho& \pho& \pho& \pho\\
			\pho& \pho& \pho& \Rt{2pt}& \Bot{2pt}\pho& \pho& \pho& \pho& \pho& \pho\\
			\pho& \pho& \pho& \pho& \Rt{2pt}\otimes& \Bot{2pt}\pho& \pho& \pho& \pho& \pho\\
			\pho& \pho& \pho& \otimes& \pho& \Rt{2pt}\pho& \Bot{2pt}\pho& \pho& \pho& \pho\\
			\red\otimes& \pho& \pho& \pho& \pho& \pho& \Rt{2pt}\pho& \Bot{2pt}\pho& \pho& \pho\\
			\pho& \red\otimes& \pho& \pho& \pho& \pho& \pho& \Rt{2pt}& \Bot{2pt}\pho& \pho\\
			\pho& \pho& \red\otimes& \pho& \pho& \pho& \pho& \pho& \red\Rt{2pt}\otimes& \pho\\
		}
	\end{equation*}
	Recall the notion of $M^D_{i,j}$ from page~\pageref{deno:M_D_ij} and the definition of $D$-regular roots, i.e., of the set $\Ro_D$ from (\ref{formula:S_D_R_D}). Recall also the notion of $\alpha_j$ and $\beta_i$ for $(i,j)\in\Phi(n)$ from (\ref{formula:alpha_j_beta_i}). In the sequel, we will write $\beta_i^1$ instead of $\beta_i$; we will also introduce the new polynomials $\beta_i^2$, $\gamma_1$ and $\gamma_2$ by putting
	\begin{equation*}
		\beta_i^2(x)=\sum_{l=3}^{i-1}x_{i,l}x_{l,2},~\gamma_1(x)=\sum\limits_{k=3}^{n-1}x_{k,1}x_{n,k},~\gamma_2(x)=\sum_{l=2}^{n-1}x_{n,l}x_{l,2}\text{ for }x\in U
	\end{equation*}

	Define also the element $x_D(\vfi)$ completely similar to (\ref{formula:x_D_vfi}):
	$$x_D(\vfi)=1_n+\sum_{(i,j)\in D}\vfi(i,j)e_{i,j}\in U,$$ where $\vfi\colon D\to\Fp_q^{\times}$ is an arbitrary map.\newpage
	
	For a map $\vfi\colon D\to\Fp_q^{\times}$, define the following subvarieties of $U(\Kp)$:
	
	\begin{equation}\label{formula:K_D_vfi}
		\begin{split}
			M_D(\vfi)(\Kp)&=\{x\in U(\Kp)\mid M_{i,j}^{D}(x)=M_{i,j}^{D}(x_D(\vfi))\text{ for } (i,j)\in R(D)\},\\
			\Ko_{D}(\vfi)(\Kp)&=
			\begin{cases}
				M_D(\vfi)(\Kp),\text{ if } D\supset D_1,\\
				\{x\in M_D(\vfi)(\Kp)\mid \beta_{i}^{2}(x)=0,~ \alpha_{j}(x)=0 \text{ for } (i,j)\in D',~\gamma_2(x)=\gamma_2(x_D(\vfi))\},\text{ if } D\supset D_2,\\
				\{x\in M_D(\vfi)(\Kp)\mid \beta_{i}^{1}(x)=0,~ \alpha_{j}(x)=0 \text{ for } (i,j)\in \overline{D},~\gamma_1(x)=\gamma_1(x_D(\vfi))\},\text{ if } D\supset D_3,\\
				\{x\in M_D(\vfi)(\Kp)\mid \beta_{i}^{2}(x)=0,~ \alpha_{j}(x)=0 \text{ for } (i,j)\in D',~\gamma_2(x)=\gamma_2(x_D(\vfi))\},\text{ if } D\supset D_4.\\
				\{x\in M_D(\vfi)(\Kp)\mid \beta_{i}^{2}(x)=0,~ \alpha_{j}(x)=0 \text{ for } (i,j)\in D',~\gamma_2(x)=\gamma_2(x_D(\vfi))\},\text{ if } D\supset D_5.
			\end{cases}
		\end{split}
	\end{equation}
	Note that $$M_{i,j}^{D}(x_D(\vfi))=\pm\prod_{(r,s)\in D^{(i,j)}\cup\{(i,j)\}}\vfi(r,s).$$ As above, we denote by $\Ko_D(\vfi)\subset U$ the set of $\Fp_q$-points of $\Ko_D(\vfi)(\Kp)$.

	Our first main result, a description of the support of $\chi$, is formulated in the following theorem\break (cf. Proposition~\ref{prop:0_1_depth_support} and Theorem~\ref{theo:0_1_depth_character} (i)).
	\mtheo{\textup{i)} For an arbitrary $2$-regular subset $D$ and an arbitrary map $\vfi\colon D\to\Fp_q^{\times}$, the conjugacy class of $x_D(\vfi)$ coincides with $\Ko_D(\vfi)$. \textup{ii)} Assume\label{theo:2_depth_support} that $x\in U$ satisfies $\chi(x)\neq0$. Then there exist the unique $2$-regular subset $D$ and the unique map $\vfi$ such that $x\in \Ko_D(\vfi)$. More precisely,
		\begin{enumerate}
			\item[\textup1.]if $x_{3,2}\neq0,x_{n,n-1}\neq0$ then $D\supset D_4$ or $D\supset D_5$, and the map $\vfi$ is so that\\ 
			$\xi_{3,n}^2\vfi(n-1,2)\vfi(n,n-1)^2-\xi_{2,n-2}\xi_{1,n-1}\vfi(n-2,1)\vfi(3,2)^2=0$\textup;
			\item[\textup2.]if $x_{3,2}=0,
			x_{n,n-1}\neq0$ then $D\supset D_3$, and the map $\vfi$ is so that\\
			$\xi_{2,n-2}\vfi(3,1)=\xi_{3,n}\vfi(n,n-1)$\textup;
			\item[\textup3.]if $x_{3,2}\neq0,x_{n,n-1}=0$ then
			$D\supset D_2$, and the map $\vfi$ is so that\\
			$\xi_{2,n-2}\vfi(3,2)=\xi_{3,n}\vfi(n,n-2)$\textup;
			\item[\textup4.]if $x_{3,2}=0,x_{n,n-1}=0$
			then $D\supset D_1$.\end{enumerate}}
	
	It remains to calculate explicitly the value of $\chi$ on all conjugacy classes described above. To formulate the answer, we need the following additional notation. Let $D$ be an arbitrary 2-regular subset, denote
	\begin{equation*}
		\begin{split}
			&\Phi_{\mathrm{\mathrm{2reg}}}=\{(i,j)\in\Phi(n)\mid i>n-i+1,(i,j)\neq (n,2)\}\cup (n-1,2),\\
			&\Phi^{'}=\{(i,j)\in\Phi_{\mathrm{\mathrm{2reg}}}\ |\ i\neq n,i\neq n-2\},\quad \Phi^{''}=\{(i,j)\in\Phi_{\mathrm{\mathrm{2reg}}}\ |\ i\neq n,i\neq n-1 \},
		\end{split}
	\end{equation*}
	\begin{equation}
		m_D=
		\begin{cases}
			|R(D)\cap \Phi_{\mathrm{\mathrm{2reg}}}|-2,&\text{if }D\supset D_1,\\
			|R(D)\cap \Phi^{'}|+n-3,&\text{if }D\supset D_2,\\
			|R(D)\cap \Phi^{''}|+n-3,&\text{if }D\supset D_3,\\
			|R(D)\cap \Phi^{'}|+n-3,&\text{if }D\supset D_4,\\
			|R(D)\cap \Phi^{'}|+n-3,&\text{if }D\supset D_5.
		\end{cases}\label{formula:m_D_depth_2}
	\end{equation}
	(Cf. the notion of $m_D$ from (\ref{formula:m_D_0_1}).)
	
	Our second main result is as follows (cf. Theorem~\ref{theo:0_1_depth_character} (ii)).
	\mtheo{Let\label{theo:2_depth_value} $\chi$ and $x\in\Ko_D(\vfi)\subset\Supp{\chi}$ be as above. Then
		\begin{equation*}\chi(x)=q^{m_D}\prod_{(i,j)\in D\cap D^2}\theta(\xi(i,j)\vfi(i,j)).
		\end{equation*}
	}
	
	In the rest of the section we sketch the proof of Theorem~\ref{theo:2_depth_support} (ii).
	Recall the semidirect decomposition $U=U'\rtimes\wt U$, where $\wt U\cong U_{n-1}(\Fp_q)$ (see Section~\ref{sect:Mackey_method}). As above, we denote by $\kappa$ the irreducible character of $U'$ of the form $$\kappa(x)=\theta(\xi(n-1,1)x_{n-1,1}),~x\in U'.$$ Recall also the little group $\wt U^{\kappa}=V\times\wt\Fp_q$, where $V=\{b\in\wt U^{\kappa}\mid b_{n,n-1}=0\}\cong U_{n-2}(\Fp_q)$ and $\wt\Fp_q=\{b\in\wt U\mid b_{i,j}=0\text{ if }(i,j)\neq(n,n-1)\}\cong\Fp_q$. Recall also the notion $$\wt\Phi=\Phi(n)\setminus(\Co_1\cup\Ro_{n-1}\cup\Co_{n-1}).$$
	
	Next, $\wt\psi$ denotes the irreducible character of depth 1 for $V\cong U_{n-2}(\Fp_q)$ corresponding to the 1-regular orbit $\Omega_{\wt D_2,\wt\xi}$, where $\wt D_2=\rho(D^2\cap\wt\Phi)$ and $\wt\xi=\restr{\xi}{D^2\cap\wt\Phi}\circ\rho^{-1}$. Then we set $$\psi(x)=\wt\psi(x_V)\theta(\xi(n,n-1)x_{n,n-1}),$$ where $\xi(n,n-1)=0$ if $(n,n-1)\notin D^2$, and $x_V$ is the projection of an element $x\in\wt U^{\kappa}$ to $V$. We need the following fact, cf. Proposition~\ref{prop:chi_0_1_regular_Mackey}. (A detailed proof is given in the next section.)
	\mprop{One\label{prop:chi_2_regular_Mackey} has $\chi=\chi_{\kappa,\psi}$.}
	
	\textsc{Sketch of the proof of Theorem~\ref{theo:2_depth_support} (ii)}. We consider the case $x_{3,2}\neq0$, $x_{n,n-1}\neq0$ here (in fact, technically it is the most complicated case). Recall that the set
	\begin{equation*}
		\Tu=\{g\in\wt U\mid g_{i,j}=0,\text{ if }i\neq n-1,~j<i\}=\exp\left(\langle e_{n-1,j},~(n-1,j)\in\Ro_{n-1}\setminus\Co_1\rangle_{\Fp_q}\right)
	\end{equation*}
	is a complete set or representatives for the left cosets $U/U'\wt U^{\kappa}$.
	
	Using Proposition~\ref{prop:chi_2_regular_Mackey}, standard formulas for induced characters and obvious matrix calculations, one obtain the equality
	\begin{equation*}
		\chi(x)=c\sum_{\{t\in\Tu,~x^t\in U'\wt U^{\kappa}\}}\kappa(x^t_{U'})\psi(x^t_{\wt U^{\kappa}}),
	\end{equation*}
	where $x^t=t^{-1}xt$\label{page:x_t}, and $c$ is a certain constant from $\Cp^{\times}$ (in fact, a power of $q$). Furthermore, $(x^t)_{i,j}=x_{i,j}$ if $i\leq n-2$, while for all $2\leq j\leq n-2$ one has
	\begin{equation}
		\begin{split}
			(x^t)_{n-1,j}&=x_{n-1,j}-\sum_{l=j+1}^{n-2}t_{n-1,l}x_{l,j},\\
			(x^t)_{n,j}&=x_{n,j}+x_{n,n-1}t_{n-1,j}.
		\end{split}\label{formula:x_t}
	\end{equation}
	
	The key observation is that $\chi(x)\neq0$ implies $\psi(x^t_{\wt U^{\kappa}})\neq0$ for at least one $t\in\Tu$. But the description of $\Supp{\wt\psi}$ from Section~\ref{sect:orbit_method} immediately shows that there exists a 1-regular subset $\wt D\subset\Phi(n-2)$ and a map $\wt\vfi\colon\wt D\to\Fp_q^{\times}$. such that $x$ satisfies all the equations for $\Ko_{\wt D}(\wt\vfi)$ (after the identification $V\cong U_{n-2}(\Fp_q)$). In particular, $$M^D_{i,j}(x)=M^D_{i,j}(x_{\wt D}(\wt\vfi))\text{ for }(i,j)\in R(D),$$ where $D=\rho^{-1}(\wt D)\setminus\Co_2$.
	
	To obtain a part of remaining equations, namely, $M^D_{n-1,j}(x)=0$ for $2\leq j\leq n-2$, note that the condition $x^t\in U'\wt U^{\kappa}$ can be simply written as $$(x^t)_{n-1,j}=0\text{ for }2\leq j\leq n-2.$$ It can be checked that the Kronecker--Capelli criterion on the existence of a solution for this system of linear equations in variables $t_{n-1,j}$ gives us the required equations $M^D_{n-1,j}(x)=0$. Using the description of $\Ko_{\wt D}(\wt\vfi)$, one can also obtain the equations $\alpha_j(x)=0$ for $(i,j)\in D$.
	
	Next, it is well known that $\sum_{c\in\Fp_q}\theta(c)=0$. Hence, if we will express some of variables $t_{n-1,j}$ as linear combinations of the other variables and substitute these expressions into $\chi(x)$, the coefficients at the remaining variables should be zero. After some delicate matrix calculations, this gives us the remaining equations for $\Ko_D(\vfi)$, i.e., $M_{i,1}(x)=0$ and $\beta_i(x)=0$. All other conditions well be proved in more delicate way, see below.\hfill{}$\square$

	\sect{Proofs: conjugacy classes and a Mackey decomposition}\label{sect:proofs_classes}
	
	This section contains detailed proofs of some of our main results stated in the previous section. Namely, we present an explicit description for the conjugacy class of $x_D(\vfi)$ (i.e., prove Theorem~\ref{theo:2_depth_support}~(i)), and find a ``Mackey decomposition'' for a character of depth 2 in the spirit of Section~\ref{sect:Mackey_method} (i.e., prove Proposition~\ref{prop:chi_2_regular_Mackey}).
	\sst{Proof of Theorem~\ref{theo:2_depth_support} (i)} Here we
	present an explicit description of the conjugacy class of an element $x_D(\varphi)$. Remind that we split an arbitrary subset $D$ into the union
	$D=\overline{D}\sqcup D^*$, where $D^*=\{(i,j)\in D\mid i=3,n-1 \text{ or } n\}$.
	
	\lemmp{For an arbitrary $2$-regular subset $D\subset\Phi(n)$ and an arbitrary map $\vfi\colon D\to\Fp_q^{\times}$\textup,\break the subvariety~$\Ko_D(\vfi)(\Kp)$ is $U(\Kp)$-invariant\textup, i.e.\textup, if $g\in\Ko_D(\vfi)(\Kp)$ then
		$xgx^{-1}\in\Ko_D(\vfi)(\Kp)$ for\break all $x\in U(\Kp)$.}{Since the set $D$ is a
		rook placement, all $M_{i,j}^{D}$ are $U(\Kp)$-invariant (the proof of this fact can be found in \cite[Lemma 2.1]{Andre01}).
		
		Let $g=(y_{i,j})\in \Ko_D(\vfi)(\Kp)$. For $\lambda\in\Kp$, denote $x_{r,s}(\lambda)=1+\lambda e_{r,s}$.
		Since every element $x\in U(\Kp)$ can be written as a product
		$x=x_{r_1,s_1}(\lambda_1)\ldots x_{r_k,s_k}(\lambda_k)$ for a certain
		$k$ and $(r_i,s_i)\in\Phi(n)$, $\lambda_i\in \Kp$, it is enough to
		prove that if $x=x_{r,s}(\lambda)$, then $xgx^{-1}\in \Ko_D(\vfi)(\Kp)$. But if
		$x$ has such a form then
		\begin{equation*} (xgx^{-1})_{i,j}=\begin{cases} y_{i,j},&\text{if
				}i\neq r\text{ and
				}j\neq s,\\
				y_{i,j},&\text{if }i=r\text{ and }j\geq s\text{, or }j=s\text{ and }i\leq r,\\
				y_{r,j}+\lambda y_{s,j},&\text{if }i=r,j<s,\\
				y_{i,s}-\lambda y_{i,r},&\text{if }j=s,i>r.
			\end{cases}
		\end{equation*}
		The proof of the invariance can be performed now by direct consideration of possible
		values of $r$ and $s$.
		
		For example, consider the polynomial $\alpha_{j}$. If $s\leq j$, then all coordinate functions involved in this polynomial are invariant themselves; this is also true, if
		$r=n$, $s\leq n-m+1$. If $r=n$, $s> n-m+1$, then
		\begin{equation*}
			\begin{split}
				\alpha_{j}(xgx^{-1})&=\sum_{l=n-m+1}^{s-1}(y_{n,l}+\lambda
				y_{s,l})y_{l,j}+
				\sum_{l=s}^{n-1}y_{n,l}y_{l,j}=\\
				&=\alpha_{j}(g)+\lambda\sum_{l=n-m+1}^{s-1}
				y_{s,l}y_{l,j}=\alpha_{j}(g),
			\end{split}
		\end{equation*} because
		$y_{s,l}=0$ for all $l\geq n-m+1\geq m$.
		
		If $r<n$, $s= j$, then $\alpha_{j}$ is invariant for the same
		reasons. Finally, if $r<n$, $s>j$, then
		\begin{equation*}\begin{split} \alpha_{j}(xgx^{-1})&= \sum_{l\neq
					r,s}y_{n,l}y_{l,j}+(y_{r,j}+\lambda
				y_{s,j})y_{n,r}+y_{s, j}(y_{n,s}-\lambda y_{n,r})=\\
				&=\alpha_{j}(g)+\lambda\cdot(y_{s,j}y_{n,r}-y_{s,
					j}y_{n,r})=\alpha_{j}(g).
		\end{split}\end{equation*}
		The invariance of $\beta_{i}^{1}$, $\beta_{i}^2$, $\gamma_1$ and $\gamma_2$ can be proved similarly.}
	\label{lemm_proofs_inv}
	
	This means that $\Ko_D(\vfi)(\Kp)$ is a union of $U(\Kp)$-conjugacy classes. It is obvious that $x_D(\varphi)$ is contained in $\Ko_D(\vfi)(\Kp)$. Thus, the $U(\Kp)$-conjugacy class of $x_D(\vfi)$ contained in $\Ko_D(\vfi)(\Kp)$.
	
	For any root $\xi=(i,j)\in\Phi(n)$, we define its \emph{level} as
	the number $u(\xi)=i-j$. Then the formula
	\begin{equation} \xi=(i_1,j_1)\lessdot\eta=(i_2,j_2)\Leftrightarrow
		\text{ either }u(\xi)<u(\eta)\text{ or }u(\xi)=u(\eta),j_1<j_2,
	\end{equation}
	defines a linear order on the set of all roots (this order does not coincide with the order $\prec$ defined in Remark~\ref{nota:basis_defi_ideal}). For each root
	$\xi=(i,j)\in\Phi(n)$, let $I_{\xi}$ be the ideal in $\Kp[U(\Kp)]$
	generated by all $y_{\eta},\eta\lessdot\xi$, and $\xi_0=(i_0,j_0)$ be the
	maximal root in $\overline{D}$, which is less than $\xi$ (if such a root exists).
	
	Now, let $J_D(\vfi)$ be the ideal of $\Kp[U(\Kp)]=\Kp[y_{i,j}]$ generated by the left-hand sides of the defining equations for $\Ko_D(\vfi)(\Kp)$, see formula~(\ref{formula:K_D_vfi}). (We denote by $y_{i,j},~1\leq j<i\leq n$, the usual coordinate functions on $U(\Kp)$.) 
	
	\lemmp{\label{lemma:simp}For an arbitrary $2$-regular subset $D$ and an arbitrary map $\vfi\colon D\to\Fp_q^{\times}$\textup, the ideal $J_D(\varphi)$ is a prime ideal of $\Kp[y_{ij}]$.}{Assume, for example, that $D\supset D_4$ and denote $J_D(\vfi)$ by $J$. (All other cases can be considered similarly.) Consider the following transformation of coordinates:
		\begin{equation}
			\begin{split}
				&\wt y_{n-1,j}=\alpha_{j},\quad
				\wt y_{i,3}=\beta_{i}^{2},\quad(i,j)\in \overline{D},\\
				&\wt y_{ij}=\Delta_{ij}-c_{ij},\quad(i,j)\in D,\\
				&\wt y_{ij}=\Delta_{ij},\quad(i,j)\in R(D).\\
				&\wt y_{n,3}=\gamma_2.
			\end{split}\label{formula:zamena}
		\end{equation}
		(it's easy to see that $J=\langle \wt y_\xi\rangle_{\xi\in B}$, where
		$B\subset\Phi(n)$ denotes the set of all roots $(i,j)$ from the
		right-hand side of (\ref{formula:zamena})). Note that for any $\xi\in B$ we
		have
		\begin{equation}\wt y_{\xi}\equiv y_{\xi}^0\cdot
			y_{\xi}+a_{\xi}\pmod{I_{\xi}},\label{treug}\end{equation} where
		$y_{\xi}^0$ is an invertible element of $\Kp[y_{ij}]/J$, and
		$a_{\xi}\in\Kp$ is a certain scalar. We can see that for this roots:
		\begin{equation*}
			\begin{split}
				&\wt y_{n-1,j}=y_{n-1,j}y_{n,n-1}+\dots,\ (i,j)\in \overline{D},\\
				&\wt y_{i,3}=y_{i,3}g_{3,2}+\dots,\ (i,j)\in \overline{D},\\
				&\wt y_{i,j}=y_{i,j}\Delta_{i_0,j_0}+\dots,\ (i,j)\in D,j\neq 1\\
				&\wt y_{i,1}=y_{i,1}g_{3,2}\Delta_{i_0,j_0}+\dots,\ (i,j)\in D,\\
				&\wt y_{i,j}=y_{i,j}\Delta_{i_0,j_0}+\dots,\ (i,j)\in R(D),j\neq 1\\
				&\wt y_{i,1}=y_{i,1}g_{3,2}\Delta_{i_0,j_0}+\dots,\ (i,j)\in R(D),\\
				&\wt y_{n,3}=y_{n,3}y_{3,2}+\dots.
			\end{split}
		\end{equation*}
		where for any $\wt y_{\xi}$, $\xi\in B$, dots denote elements equal
		to zero modulo $I_{\xi}$ and scalars (we assume that
		$\Delta_{i_0,j_0}=1$ if the root $\xi_0$ does not exist for a given
		$\xi\in\Phi(n)$). But one can easily obtain the following equalities
		modulo $J$:
		\begin{equation*}\begin{split} 
				&y_{n,n-1}\equiv c_{n,n-1}\neq 0,\quad y_{3,2}\equiv c_{3,2}\neq 0\\
				&\Delta_{i_0,j_0}\equiv c_{i_0,j_0}\neq 0,\ (i_0,j_0)\in \overline{D}.
			\end{split}
		\end{equation*}
		This concludes the proof of
		(\ref{treug}).
		
		Hence, in the quotient algebra $\Kp[y_{ij}]/J$, all $y_{\xi}$ are polynomial in $\wt
		y_{\xi}$ (for each $\xi\in B$). Consequently,
		$$\Kp[y_{ij}]/J=\left.\Kp[y_{\xi}]_{\xi\in\Phi(n)}
		\right/\langle\wt y_{\xi}\rangle_{\xi\in B} \cong\left.\Kp[\{y_{\xi}\}_{\xi\notin B}\cup\{\wt y_{\xi}\}_{\xi\in
			B}]\right./\langle\wt y_{\xi}\rangle_{\xi\in B}\cong \Kp[\wt
		y_{\xi}]_{\xi\notin B}.$$ In particular, $\Kp[y_{ij}]/J$ is a
		domain, so $J$ is a prime ideal.}\label{lemm_proofs_prim}
	
	This Lemma shows that the subvariety $\Ko_D(\vfi)$ of $U(\Kp)$ corresponding to the ideal $J_D(\vfi)$ is irreducible.
	\lemmp{Let\label{lemm:dim_stab}
		$\mathcal C=\{g\in U(\Kp)\mid
		gx_D(\varphi)=x_D(\varphi)g\}$ be the stabilizer
		of $x_D(\varphi)$ in $U(\Kp)$. Then $\dim \mathcal C=\codim\Ko_D(\vfi)(\Kp)$.}{First, we will present the defining equations
		for $\mathcal C$. We~put $$D^{+}=\{(i,j)\in D\mid i<n-2,~j>3\}$$ and, for any $\xi=(i,j)\in\Phi(n)$,
		\begin{equation*}
			\begin{split}
				&\Phi_{\xi}=\{(i^{'},j^{'})\mid i^{'}=j,~1\leq j^{'}<j\},\quad \Phi^{\xi}=\{(i^{'},j^{'})\mid j^{'}=i,~j<i^{'}\leq n\}.
			\end{split}
		\end{equation*}
		We will consider all possible cases for $D\supset D_i$, $1\leq i\leq 5$. 
		
		First, let the subset $D$ contain $D_1$. It is easy to check that $\mathcal
		C$ is defined by the following equations:
		\begin{equation*}\begin{split}
				&a_{n-2,2}y_{2,1}=0,~a_{n-2,2}y_{n-1,n-2}=0,\\
				&a_{n-2,2}y_{2,1}+a_{n,3}y_{3,1}-a_{n-1,1}y_{n,n-1}=0,\\
				&a_{n,3}y_{3,2}-a_{n-2,2}y_{n,n-2}=0,\\
				&y_{i,j}=0,~(i,j)\in\bigsqcup_{\xi\in D^{+}}(\Phi_{\xi}\cup\Phi^{\xi})\\
		\end{split}\end{equation*}
		(here and further we will write $a_{i,j}=\varphi(i,j)\in \Fp_q^{\times}$ if $(i,j)\in D$, else $a_{i,j}=0$). Hence, the number of the equations on $\mathcal{C}$ is equal to $|\Phi(n)\backslash R(D)|$, and all these linear equations are linearly independent.

		Second, let the subset $D$ contain $D_2$. Put 
		\begin{equation*}\begin{split}
				\Phi_{\alpha}&=\{(n-2,j)\mid (n-j+1,j)\notin D^{+},~(j,n-j+1)\notin D^{+},~3<j<n-2\},\\
				\Phi_{\beta}&=\{(i,3)\mid (i,n-i+1)\notin D^{+},~(n-i+1,i)\notin D^{+},~3<i<n\}.  
			\end{split}    
		\end{equation*}
		One can easily see that $\mathcal C$ is given by the following equations:
		\begin{equation*}\begin{split}
				&y_{i,j}=0,~(i,j)\in\bigsqcup_{\xi\in D^{+}}(\Phi_{\xi}\cup\Phi^{\xi}), i\neq n, j\neq 2,\\
				&y_{i,3}=0,~ (i,3)\in\Phi_{\beta},\\
				&y_{n-2,j}=0,~ (n-2,j)\in\Phi_{\alpha},\\
				&a_{i,j}y_{j,2}-y_{i,3}a_{3,2}=0,~ (i,j)\in D^{+},\\
				&a_{i,j}y_{n,j}-y_{n-2,j}a_{n,n-2}=0,~ (i,j)\in D^{+},\\
				&a_{n-1,1}y_{n,n-1}-a_{n,n-2}y_{n-2,1}=0,~ y_{n,3}a_{3,2}-y_{n-2,2}a_{n,n-1}=0.
		\end{split}\end{equation*}
		The number of the equations in this case is equal to $|\Phi(n)\backslash R(D)|-2|D^{+}|-1$, and all these linear equations are linearly independent. 
		
		Next, assume that the subset $D$ contains $D_3$. Put 
		\begin{equation*}\begin{split}
				\Phi_{\alpha}&=\{(n-1,j)\mid (j,n-j+1)\notin \overline{D},~(n-j+1,j)\notin D^{+},~3<j<n-1\},\\
				\Phi_{\beta}&=\{(i,3)\mid (n-i+1,i)\notin D^{+},~(i,n-i+1)\notin \overline{D},~3<i<n-1\}.  
			\end{split}
		\end{equation*}
		The stabilizer $\mathcal C$ is described by the following equations:
		\begin{equation*}\begin{split}
				&y_{i,j}=0,~(i,j)\in\bigsqcup_{\xi\in D^{*}}(\Phi_{\xi}\cup\Phi^{\xi}), i\neq n, j\neq 1,\\
				&y_{i,3}=0,~ (i,3)\in\Phi_{\beta},\\
				&y_{n-1,j}=0,~ (n-1,j)\in\Phi_{\alpha},\\
				&a_{i,j}y_{j,1}-y_{i,3}a_{3,1}=0,~ (i,j)\in \overline{D},\\
				&a_{i,j}y_{n,j}-y_{n-1,j}a_{n,n-1}=0,~ (i,j)\in \overline{D},\\
				&y_{n-1,3}=0,~ y_{n,3}a_{3,1}-y_{n-1,1}a_{n,n-1}=0.
		\end{split}\end{equation*}
		The number of the equations in this case is equal to $|\Phi(n)\backslash R(D)|-|D^{+}|-|\overline{D}|-1$, and all these linear equations are linearly independent. 
		
		Now, assume that the subset $D$ contains $D_4$. Put
		\begin{equation*}\begin{split}
				\Phi_{\alpha}&=\{(n-1,j)\mid (j,n-j+1)\notin D^+,~(n-j+1,j)\notin D^{+},~3<j<n-1\},\\
				\Phi_{\beta}&=\{(i,3)\mid (n-i+1,i)\notin D^{+},~(i,n-i+1)\notin D^+,~3<i<n-1\}.  
			\end{split}
		\end{equation*}
		It is easy to check that $\mathcal C$ is defined by the following equations: 
		\begin{equation*}\begin{split}
				&y_{i,j}=0,~(i,j)\in\bigsqcup_{\xi\in D^{*}}(\Phi_{\xi}\cup\Phi^{\xi}), i\neq n, j\neq 2,\\
				&y_{i,3}=0,~ (i,3)\in\Phi_{\beta},\\
				&y_{n-1,j}=0,~ (n-1,j)\in\Phi_{\alpha},\\
				&a_{i,j}y_{j,2}-y_{i,3}a_{3,2}=0,~ (i,j)\in D^+,\\
				&a_{i,j}y_{n,i}-y_{n-1,j}a_{n,n-1}=0,~ (i,j)\in D^+,\\
				&y_{n-1,3}=0,~ a_{n,3}y_{3,1}-y_{n-1,1}a_{n-1,1}+y_{n-1,1}a_{n,n-1}=0,\\
				&a_{n,3}y_{3,2}+a_{n,n-1}y_{n-1,2}-a_{3,2}y_{n,3}=0.
		\end{split}\end{equation*}
		The number of the equations in this case is equal to $|\Phi(n)\backslash R(D)|-2|D^{+}|-1$, and all these linear equations are linearly independent.
		
		Finally, assume that the subset $D$ contains $D_5$. Put
		\begin{equation*}\begin{split}
				\Phi_{\alpha}&=\{(n-1,j)\mid (j,n-j+1)\notin D^+,~(n-j+1,j)\notin D^{+},~3<j<n-1\},\\
				\Phi_{\beta}&=\{(i,3)\mid (n-i+1,i)\notin D^{+},~(i,n-i+1)\notin D^+,~3<i<n-1\}.  
			\end{split}
		\end{equation*}
		It is easy to check that $\mathcal C$ is defined by the following equations: 
		\begin{equation*}\begin{split}
				&y_{i,j}=0,~(i,j)\in\bigsqcup_{\xi\in D^{*}}(\Phi_{\xi}\cup\Phi^{\xi}), i\neq n, j\neq 2,\\
				&y_{i,3}=0,~ (i,3)\in\Phi_{\beta},\\
				&y_{n-1,j}=0,~ (n-1,j)\in\Phi_{\alpha},\\
				&a_{i,j}y_{j,2}-y_{i,3}a_{3,2}=0,~ (i,j)\in D^+,\\
				&a_{i,j}y_{n,i}-y_{n-1,j}a_{n,n-1}=0,~ (i,j)\in D^+,\\
				&y_{n-1,3}=0,~ a_{n,3}y_{3,1}-y_{n,n-2}a_{n-2,1}+y_{n-1,1}a_{n,n-1}=0,\\
				&a_{n,3}y_{3,1}+a_{n,n-1}y_{n-1,1}-a_{n-2,1}y_{n,n-2}=0,~y_{2,1}a_{n-1,2}-y_{n-1,n-2}a_{n-2,1}.
		\end{split}\end{equation*}
		The number of the equations in this case is equal to $|\Phi(n)\backslash R(D)|-2|D^{+}|-2$, and all these linear equations are linearly independent.
		
		It is obvious that the dimension of $\Ko_D(\vfi)(\Kp)$ (one can compute this dimension using Lemma~\ref{lemma:simp}) is equal to the number of the equations in all cases. This concludes the proof.}
	
	Clearly, the $U_n(\Kp)$-conjugacy class of $x_D(\vfi)$ is an irreducible subvariety of $U_n(\Kp)$ contained in the irreducible $U_n(\Kp)$-invariant subvariety $\Ko_D(\vfi)(\Kp)$, because $x_D(\vfi)\in\Ko_D(\vfi)(\Kp)$. The dimensions of these subvarieties are equal by Lemma~\ref{lemm:dim_stab}, so the subvarieties coincide. It follows from Lemma~\ref{lemma:simp} that $\Ko_D(\vfi)(\Kp)$ as an affine variety is defined over $\Fp_q$ and $\Fp_q$-isomorphic to the affine space of dimension $\dim(\Ko_D(\vfi)(\Kp))$. Hence, $\Kp_D(\vfi)$, being the set of $\Fp_q$-points of $\Ko_D(\vfi)(\Kp)$, consists of $q^{\dim(\Ko_D(\vfi)(\Kp))}$ elements. On the other hand, the proof of Lemma~\ref{lemm:dim_stab} in fact shows that the $U_n$-conjugacy class of $x_D(\vfi)$ contains $q^{\codim\Co}$ elements. But $x_D(\vfi)\in\Ko_D(\vfi)$ and $\dim\Ko_D(\vfi)=\codim\Co$, thus, the conjugacy class of $x_D(\vfi)$ coincides with $\Ko_D(\vfi)$. The proof of Theorem~\ref{theo:2_depth_support} (i) is complete. 
	

	\sst{Proof of Proposition~\ref{prop:chi_2_regular_Mackey}} We need to introduce some notation very similar to the notation used in Section~\ref{sect:Mackey_method}. As above, let $\Omega_{D^2,\xi}=\Omega_{f_{D^2,\xi}}$ be a 2-regular orbit and $\chi=\chi_{D^2,\xi}$ be the corresponding character of depth~2. Recall the notion of $U'$ and $\wt U$ from Section~\ref{sect:Mackey_method}, so that $U=U'\rtimes\wt U$. Recall the notion of $\kappa$, $\wt U^{\kappa}$, $V$, $\wt\Fp_q$, $\Tu$, $\wt\Phi$, $\wt D_2$, $\wt\xi$, $\rho$, $\wt\psi$ and $\psi$ from the previous section, as well as the notion $\vt$ for the Lie algebra of $V$. Slightly abusing the notation, we will denote all natural isomorphisms $V\cong U_{n-2}(\Fp_q)$, $\vt\cong\ut_{n-2}(\Fp_q)$ and $\vt^*\cong\ut_{n-2}(\Fp_q)^*$ by the same letter $\rho$. Put also $$\wt f=\sum_{(i,j)\in \wt D_2\cap\wt\Phi}\wt\xi(i,j)e_{j,i}\in\vt^*,$$ so that $\Omega_{\wt D,\wt\xi}$ is the coadjoint $U_{n-2}(\Fp_q)$-orbit of the linear form $\rho(\wt f)\in\ut_{n-2}(\Fp_q)^*$.
	
	Now, recall from Lemma~\ref{lemm:polarization} that $\pt=\langle e_{i,j},~(i,j)\in\Phi(n)\setminus\Mo\rangle_{\Fp_q}\subset\ut$ is a polarization for $f=f_{D^2,\xi}$. Note also that $$\wt\pt=\langle e_{i,j},~(i,j)\in\Phi(n)\setminus(\Mo\cup\Co_1\cup\Ro_{n-1}\cup\{(n,n-1)\})\rangle_{\Fp_q}\subset\vt$$ is a polarization for $\wt f$, see Example~\ref{exam:0_1_regular_orbits} (ii). On the picture below we drew schematically the situation for $n=10$. We keep the notations from Example~\ref{exam:2_regular_orbits}. As in Example~\ref{exam:0_1_decomposition_Mackey}, we marked roots from $\Co_1$ gray, roots from $\Ro_{n-1}$ yellow, and the root $\{(n,n-1)\}$ green, so that the roots from $\wt\Phi$ are exactly the white roots. Since the roots from $\Mo$ are marked by minuses, we hope that it is clear that $\pt\cap\vt=\wt\pt$.
	\begin{equation*}
		\mymatrix{
			\Bot{2pt}\pho& \pho& \pho& \pho& \pho& \pho& \pho& \pho& \pho& \pho\\
			\gray\Rt{2pt}+& \Bot{2pt}\pho& \pho& \pho& \pho& \pho& \pho& \pho& \pho& \pho\\
			\gray+& \Rt{2pt}+& \Bot{2pt}\pho& \pho& \pho& \pho& \pho& \pho& \pho& \pho\\
			\gray+& +& \Rt{2pt}+& \Bot{2pt}\pho& \pho& \pho& \pho& \pho& \pho& \pho\\
			\gray+& +& +& \Rt{2pt}+& \Bot{2pt}\pho& \pho& \pho& \pho& \pho& \pho\\
			\gray+& +& +& +& \Rt{2pt}\otimes& \Bot{2pt}\pho& \pho& \pho& \pho& \pho\\
			\gray+& +& +& \otimes& -& \Rt{2pt}-& \Bot{2pt}\pho& \pho& \pho& \pho\\
			\gray+& \otimes& -& -& -& -& \Rt{2pt}-& \Bot{2pt}\pho& \pho& \pho\\
			\gray\otimes& \yellow-& \yellow-& \yellow-& \yellow-& \yellow-& \yellow-& \Rt{2pt}\yellow-& \Bot{2pt}\pho& \pho\\
			\gray\bullet& \bullet& \otimes& -& -& -& -& \square& \Rt{2pt}\green\square& \pho\\}
	\end{equation*}
	
	For brevity, denote $A=U'$, $B=\wt U$, so that $\wt U^{\kappa}=B^{\kappa}=V\times\wt\Fp_q$. Note that, given $x\in AB^{\kappa}$,
	\begin{equation*}
		\kappa_0(x)=\kappa(x_A)=\theta(\xi(n-1,1)x_{n-1,1})=\theta(f(\ln(\pi^{AB^{\kappa}}_A(x)))).
	\end{equation*}
	At the same time, by definition of the character associated with a coadjoint orbit (see Section~\ref{sect:orbit_method}),
	\begin{equation*}
		\begin{split}
			\psi_0(x)=\psi(x_{B^{\kappa}})=\theta(\xi(n,n-1)x_{n,n-1})\wt\psi(x_V)=\theta(\xi(n,n-1)x_{n,n-1})\left(\Ind{\wt P}{V}{(\theta\circ\wt f\circ\ln)}\right)(x_V),
		\end{split}
	\end{equation*}
	where we put $\wt P=\exp(\wt\pt)$ and $x_V=\pi^{AB^{\kappa}}_V(x)=\pi^{B^{\kappa}}_V(\pi^{AB^{\kappa}}_{B^{\kappa}}(x))$. Since $AB^{\kappa}=(A\rtimes V)\rtimes\wt\Fp_q$, we can also write $\pi^{AB^{\kappa}}_V=\pi^{B^{\kappa}}_V\circ\pi^{AB^{\kappa}}_{B^{\kappa}}=\pi^{AV}_{V}\circ\pi^{AB^{\kappa}}_{AV}$.
	
	According to \cite[Proposition 1.2]{Lehrer74}, given a finite group $A\rtimes B$, a subgroup $C\subset B$, and a character $\lambda$ of $C$, one has $(\Ind{C}{B}{\lambda})\circ\pi^{AB}_B=\Ind{AC}{AB}{(\lambda\circ\pi^{AC}_C)}$. Hence,
	\begin{equation*}
		\begin{split}
			(\Ind{\wt P}{V}{(\theta\circ\wt f\circ\ln)})\circ\pi^{AV}_V=\Ind{A\wt P}{AV}{(\theta\circ\wt f\circ\ln\circ\pi^{A\wt P}_{\wt P})}.
		\end{split}
	\end{equation*}
	Denote this character of $AV$ by $\zeta$. Since $$M=\exp\langle e_{i,j},~(i,j)\in\Mo\rangle_{\Fp_q}=\{g\in U\mid g_{r,s}=0\text{ for }(r,s)\notin\Mo\}$$ is a complete system of representatives for the left cosets $AV/A\wt P$, the usual formula for induced character shows that, for $g\in AV$,
	\begin{equation*}
		\begin{split}
			\zeta(g)=\sum_{\{h\in M,~h^{-1}gh\in A\wt P\}}\theta(\wt f(\ln(\pi^{A\wt P}_{\wt P}(h^{-1}gh)))),
		\end{split}
	\end{equation*}
	and one can write $f$ instead of $\wt f$ in this formula for an obvious reason.
	
	Next, we note that $(aba^{-1})_{n,n-1}=b_{n,n-1}$ for all $a,~b\in U$, and denote for simplicity $x_{AV}=\pi^{AB^{\kappa}}_{AV}(x)$, $x_{\Fp_q}=x_{n,n-1}=\pi^{AB^{\kappa}}_{\wt\Fp_q}(x)$, so that $\xi(n,n-1)x_{n,n-1}=f(\ln(x_{\Fp_q}))$. Clearly, for a finite group $G=A\rtimes B$ one has $b^{-1}g_Ab=(b^{-1}gb)_A$ for all $g\in G$, $b\in B$. We obtain
	\begin{equation*}
		\begin{split}
			\psi_0(x)&=\theta(f(\ln(x_{\Fp_q})))\zeta(x_{AV})=\sum_{\{h\in M,~h^{-1}x_{AV}h\in A\wt P\}}\theta(f(\ln(\pi^{A\wt P}_{\wt P}(h^{-1}x_{AV}h))+\ln(x_{\Fp_q})))\\
			&=\sum_{\{h\in M,~(h^{-1}xh)_{AV}\in A\wt P\}}\theta(f(\ln(\pi^{A\wt P}_{\wt P}((h^{-1}xh)_{AV})+\ln((h^{-1}xh)_{\Fp_q}))).
		\end{split}
	\end{equation*}
	
	Since $AB^{\kappa}=(A\rtimes V)\rtimes\wt\Fp_q$, we have $P=A\wt P\wt\Fp_q=A\rtimes(\wt P\times\wt\Fp_q)=(A\rtimes\wt P)\rtimes\wt\Fp_q$, so one can write $\pi^P_{\wt P}=\pi^{A\wt P}_{\wt P}\circ\pi^P_{A\wt P}=\pi^{\wt P\wt\Fp_q}_{\wt P}\circ\pi^P_{\wt P\wt\Fp_q}$. Furthermore, for $y\in AB^{\kappa}$, the conditions $y_{AV}\in A\wt P$ and $y\in P$ are in fact equivalent, so we can rewrite $\psi_0(x)$ as follows:
	\begin{equation*}
		\psi_0(x)=\sum_{\{h\in M,~h^{-1}xh\in P\}}\theta\left(f\left(\ln(\pi^P_{\wt P}(h^{-1}xh))+\ln((h^{-1}xh)_{\wt\Fp_q})\right)\right).
	\end{equation*}
	
	According to the well-known Baker--Campbell--Hausdorff formula, if elements $a$, $b$ belong to a subalgebra $\at\subset\ut$, then $$\exp(a)\exp(b)=\exp(a+b+c)$$ for a certain element $c\in[\at,\at]$. Let $\at$ be the Lie algebra of $A\wt P\wt\Fp_q$. Applying the formula to the elements $a=\ln(\pi^P_{\wt P}((h^{-1}xh))$, $b=\ln((h^{-1}xh)_{\wt\Fp_q})$, we obtain the existence of $c\in[\at,\at]$ such that
	$$
	\ln(\pi^P_{\wt P}((h^{-1}xh))+\ln((h^{-1}xh)_{\wt\Fp_q})=\ln(\pi^P_{\wt P\wt\Fp_q}(h^{-1}xh))-c.
	$$
	But $P=\exp(\pt)$, so $\at$ is contained in the polarization $\pt$ for $f$, hence $f(c)=0$. We conclude that
	$$
	\psi_0(x)=\sum_{\{h\in M,~h^{-1}xh\in P\}}\theta(f(\ln(\pi^P_{\wt P\wt\Fp_q}(h^{-1}xh)))).
	$$

	As we mentioned above for a finite group $G$, $G=A\rtimes B$, one has $b^{-1}g_Ab=(b^{-1}gb)_A$ for all $b\in B$, $g\in G$. Recall also that $M$ is contained in the little group $B^{\kappa}$ for $\kappa$. Hence, for elements $x\in AB^{\kappa}$, $h\in M$, we see that
	$$
	\kappa_0(x)=\kappa(x_A)=\kappa(h^{-1}x_Ah)=\kappa((h^{-1}xh)_A)=\theta(f(\ln(h^{-1}xh))).
	$$
	Thus, using again the Baker--Campbell--Hausdorff formula and the fact that $\pt$ is a polarization for $f$, we take, for $x\in AB^{\kappa}$,
	\begin{equation*}
		\begin{split}
			\kappa_0(x)\psi_0(x)&=\sum_{\{h\in M,~h^{-1}xh\in P\}}\theta\left(f\left(\ln(\pi^P_A(h^{-1}xh))+\ln(\pi^P_{\wt P\wt\Fp_q}(h^{-1}xh))\right)\right)\\
			&=\sum_{\{h\in M,~h^{-1}xh\in P\}}\theta(f(\ln(h^{-1}xh)))=\Ind{P}{AB^{\kappa}}{(\theta\circ f\circ\ln)},
		\end{split}
	\end{equation*}
	because $M$ is a complete system of representatives for the left cosets $AB^{\kappa}/P$.
	
	Finally, $$\chi_{\kappa,\psi}=\Ind{AB^{\kappa}}{U}(\kappa_0\psi_0)=\Ind{AB^{\kappa}}{U}{(\Ind{P}{AB^{\kappa}}{(\theta\circ f\circ\ln)})}=\Ind{P}{U}{(\theta\circ f\circ\ln)}=\chi$$ by definition of $\chi$. The proof is complete.

	\sect{Proofs: the support and the value}\label{sect:support}
	
	\sst{Proof of Theorem~\ref{theo:2_depth_support} (ii)} Let $x=(x_{i,j})$ belong to the support of the character $\chi=\chi_{\Omega}$ of depth 2 corresponding to a 2-regular orbit $\Omega_{D^2,\xi}$, i.e. $\chi(x)\neq 0$. For brevity, we will right $\xi_{i,j}=\xi(i,j)$ for $(i,j)\in D^2$. According to Proposition~\ref{prop:chi_2_regular_Mackey} and usual formulas for induced characters,
	\begin{equation}
		\begin{split}
			\chi(x)&=q^c\sum\limits_{\{t\in\Tu,~txt^{-1}\in AB^{\kappa}\}}\theta(x^t_{n-1,1}\xi_{1,n-1})\psi(\wt x)\\
			&=q^c\sum\limits_{\{t\in\Tu,~txt^{-1}\in AB^{\kappa}\}}\theta(\xi_{1,n-1}(x_{n-1,1}-\sum\limits_{k=2}^{n-2}t_{n-1,k}x_{k,1}))\psi(\wt x)
			\label{formula:char}
		\end{split}
	\end{equation}
	for a certain integral constant $c$. Recall that formula~(\ref{formula:x_t}) implies that $(x^t)_{i,j}=x_{i,j}$ if $i\leq n-2$, while for all $2\leq j\leq n-2$ one has
	\begin{equation}
		\begin{split}
			(x^t)_{n-1,j}&=x_{n-1,j}-\sum_{l=j+1}^{n-2}t_{n-1,l}x_{l,j},\\
			(x^t)_{n,j}&=x_{n,j}+x_{n,n-1}t_{n-1,j}.
			\label{formula:x_t.2}
		\end{split}
	\end{equation}
	
	Since $\wt x\in\Supp{\wt\psi}$, we can use the description of the support of the 1-regular character $\wt\psi$, i.e., there exist a 1-regular subset $\wt D\subset\Phi(n-2)$ and a map $\wt\vfi\colon\wt D\to\Fp_{q}^{\times}$ such that $\wt x$ satisfies the equations for the conjugacy class of $x_{\wt D}(\wt\vfi)$. For a subset $D\in\Phi(n)$, we denote by $\Co_D$ (respectively, by $\Ro_D$) the set of the columns (respectively, of the rows) of all $(i,j)\in D$ such that $i\neq n,n-1$.
	
	First, consider the case when $x_{3,2}=x_{n,n-1}=0$. Using (\ref{formula:x_t.2}), we obtain $\wt x=x_{\wt U^{\kappa}}$, hence
	\begin{equation}
		M^D_{i,j}(x)=M^D_{i,j}(x_{D}(\vfi))\text{ for }(i,j)\in R(D)\setminus(\Co_1\cup\Ro_{n-1}),
		\label{equations:m_i,j.1}
	\end{equation}
	where $D=\rho^{-1}(\wt D)$ and $\vfi(i,j)=\wt\vfi(\rho(i,j))$ (note that $(2,1),(n,n-1)$ are not contained in $\wt D$ because $x_{3,2}=0$).
	
	Put $\vartheta=\text{max}_{(i,j)\in \overline{D}}j$ and $\ell=\text{min}_{(i,j)\in D^{+}}j$. It is easy to see from formula (\ref{formula:x_t.2}) that $(x^t)_{n-1,j}=0$ if $j>\vartheta$. We can rewrite the requirement that $x^t\in AB^{\kappa}$ as follows:
	\begin{equation}
		x_{n-1,j}-\sum_{s=j+1}^{n-2}t_{n-1,s}x_{s,j}=0\text{ for }1<j\leq\vartheta.
		\label{syst:abkappa}
	\end{equation}
	Since $\chi(x)\neq 0$, at least one summand in (\ref{formula:char}) is nonzero. In the other words, system (\ref{syst:abkappa}) should have at least one solution. Using the Kronecker--Capelli criterion, we conclude that $M_{n-1,j}^{D}(x)=0$  for $1<j\leq\vartheta$, $j\notin\Co_D.$
	
	Next, given subsets $X=\{i_1,\ldots,i_s\}$ and $Y=\{j_1,\ldots,j_s\}$ of $\Zp_{>0}$, denote by $M_X^Y(a)$ (cf. the definition of $\Delta_D^{r,s}(\lambda)$ on page~\pageref{page:delta}) the minor of a matrix $a$ with the set of columns $\sigma(j_1),\ldots,\sigma(j_s)$ and the set of rows $\tau(i_1),\ldots,\tau(i_s)$, where $\tau=\tau_Y$ (respectively, $\sigma=\sigma_X$) is the unique permutation such that $\tau(j_1)<\ldots<\tau(j_s)$ (respectively, $\sigma(i_1)<\ldots<\sigma(i_s)$). Solving system (\ref{syst:abkappa}), we see that
	\begin{equation*}
		t_{n-1,j}=\frac{\eth-\sum\limits_{\wp<k<n-1,~k\notin\Ro_D}\text{sgn}(\tau_{X_k})t_{n-1,k}M^{\Co_D}_{X_k}(x)}{M_{\Ro_D}^{\Co_D}(x)}\text{ for }j\in\Ro_D. 
	\end{equation*}
	Here $\wp=\text{min}_{(b,s)\in D}b$, $X_k$ is obtained from $\Ro_D$ by replacing $j$ by $k$, and $\eth$ is a certain constant. After substitution of these expressions into formula (\ref{formula:char}), the coefficients at the remaining variables should become zero, because $\sum_{c\in\Fp_q}\theta(c)=0$. Thus, we obtain the following equations on $x$:
	\begin{equation*}
		\begin{split}
			&x_{i,1}M_{\Ro_D}^{\Co_D}(x)-\sum\limits_{k\in\Ro_D}\text{sgn}(\tau_{X^k})x_{k,1}M_{X^k}^{\Co_D}(x)=0,~1<i<n-1,i\notin\Ro_D,\\
			&x_{i,1}=0,~2\leq i<\wp.
		\end{split}
	\end{equation*}
	Here, $X^k$ is obtained from $\Ro_D$ by replacing $k$ by $i$. These equations can be rewritten as follows:
	\begin{equation}
		M_{i,1}^{D}(x)=0,~1<i<n-1,i\notin\Ro_D.
		\label{equations:m_i,j.2}
	\end{equation}
	
	Hence, in this case, $\Supp{\chi}$ is contained in the union of the subvarieties $P_D(\vfi)$ defined by equations (\ref{equations:m_i,j.1}) and (\ref{equations:m_i,j.2}). Such a subvariety is not a single conjugacy class, but we can split it into the union of conjugacy classes. Specifically, it is easy to see that $P_D(\vfi)=\bigsqcup\limits_{\Do,\vfi}\Ko_\Do(\vfi)$, where the set $\Do\in\Phi(n)$ obtained from the set $D$ by adding one of the elements $(n,1)$ or $(n-1,1)$ or nothing, so that $\Do\supset D_1$. The map~$\vfi$ is extended to a map from $\Do$ to $\Fp_q^{\times}$ in an arbitrary way. Thus, to conclude the proof in the case $x_{3,2}=x_{n,n-1}=0$, it remains to check that $\chi(x)\neq0$ for all $x\in\Ko_{\Do}(\vfi)$. This will be done in the next subsection. 
	
	Second, consider the case when $x_{3,2}\neq 0$ and $x_{n,n-1}=0$. We can repeat arguments from the previous case, because $\wt x=x_{\wt U^{\kappa}}$. Hence, 
	\begin{equation}
		\begin{split}
			&M^D_{i,j}(x)=M^D_{i,j}(x_D(\vfi))\text{ for }(i,j)\in R(D)\setminus\Co_1,\\
			&\beta_i^2(x)=\alpha_j(x)=0\text{ for }(i,j)\in D^+,\\
			&\xi_{2,n-2}x_{3,2}=\xi_{3,n}x_{n,n-2}.
			\label{equations:2case_1part}
		\end{split}
	\end{equation}
	where $D=\rho^{-1}(\wt D)$ and $\vfi(i,j)=\wt\vfi(\rho(i,j))$ (the subset $\wt D$ is a union of the 0-regular subset $\wt D^+$ and the two-element subset $\{(2,1),~(n,n-1)\}$ in this case).
	
	As in the previous case, we need to solve system~(\ref{syst:abkappa}). The solution is as follows:
	\begin{equation*}\begin{split}
			&t_{n-1,j}=\frac{\eth-\sum\limits_{\wp'<k<n-1,~k\notin\Ro_D}\text{sgn}(\tau_{X_k})t_{n-1,k}M^{\Co_D\setminus\{2\}}_{X_k}(x)}{M_{\Ro_D\setminus \{3\}}^{\Co_D\setminus \{2\}}(x)}\text{ for }j\in\Ro_D\setminus \{3\},\\
			&t_{n-1,3}=\frac{\eth-\sum\limits_{\wp'<k<n-1,~k\notin\Ro_D}\text{sgn}(\tau_{X_k})t_{n-1,k}M^{\Co_D}_{W_k}(x)}{M_{\Ro_D}^{\Co_D}(x)}-\frac{\sum\limits_{3<j<\wp'}t_{n-1,j}x_{j,2}}{x_{3,2}}.
		\end{split}    
	\end{equation*}
	Here, $\wp'=\text{min}_{(b,s)\in D^+}b$, $X_k$ is obtained from $\Ro_D\setminus\{3\}$ by replacing $j$ by $k$, $W_k$ is obtained from $\Ro_D\setminus \{3\}$ by adding $k$, and $\eth$ is a certain constant. Arguing similarly to the previous case, we deduce that the coefficients at the remaining variables after substitutions of all expressions into formula (\ref{formula:char}) should be be zero. Thus, we get that
	\begin{equation*}
		\begin{split}
			&x_{2,1}=0,~-x_{i,1}x_{3,2}+x_{3,1}x_{i,2}=0,\text{ for }3<i<\wp',\\
			&x_{i,1}M_{\Ro_D\setminus \{3\}}^{\Co_D\setminus \{2\}}(x)-\sum\limits_{k\in\Ro_D}\text{sgn}(\tau_{X^k})x_{k,1}M_{X^k}^{\Co_D\setminus \{2\}}(x)-\frac{\text{sgn}(\tau_{W_i})x_{3,1}M_{W_i}^{\Co_D}(x)}{x_{3,2}}=0,~\wp'<i<n-1,i\notin\Ro_D,
		\end{split}
	\end{equation*}
	
	These equations can be rewritten as follows:
	\begin{equation}
		M_{i,1}^{D}(x)=0,~1<i<n-1,i\notin\Ro_D.
		\label{equations:2case_2part}
	\end{equation}
	Hence, in this case, $\Supp{\chi}$ is contained in the union of the subvarieties $P_D(\vfi)$ defined by equations (\ref{equations:2case_1part}) and (\ref{equations:2case_2part}). The subvariety $P_D(\vfi)$ is not a conjugacy class in this case, too, but it is easy to check that $P_D(\vfi)=\bigsqcup\limits_{\Do,\vfi}\Ko_{\Do}(\vfi)$, where the set $\Do\in\Phi(n)$ is obtained from the set $D$ by adding the element $(n-1,1)$ or nothing, so that $\Do\supset D_2$. 
	
	The map~$\vfi$ is extended to a map $\Do\to\Fp_q^{\times}$ in an arbitrary way so that $\xi_{2,n-2}\vfi(3,2)=\xi_{3,n}\vfi(n,n-2)$. (The latter condition comes from the fact that $\wt x$ belong to the support of the character $\psi$, see Theorem~\ref{theo:0_1_depth_character}.) Thus, to conclude the proof in the case $x_{3,2}\neq 0 $ and $x_{n,n-1}=0$, it remains to check that $\chi(x)\neq0$ for all $x\in\Ko_{\Do}(\vfi)$. This will be done in the next subsection.
	
	Now, assume that $x_{3,2}=0$ and $x_{n,n-1}\neq 0$. Since $\wt x\in\Supp{\psi}$, we see that
	\begin{equation*}
		M_{i,j}^{\wt D}(\wt x)=M_{i,j}^{\wt D}(x_{\wt D}(\wt\vfi))\text{ for }(i,j)\in R(\wt D). 
	\end{equation*}
	But $\wt x_{i,j}=x_{i,j}$ when $i<n-1$ for all $t\in\Tu$. It follows that
	\begin{equation}
		M_{i,j}^{D}(x)=M_{i,j}^{D}(x_D(\vfi))\text{ for }(i,j)\in R(D)\setminus(\Co_1\cup\Ro_n), 
		\label{equations:3case_1part}
	\end{equation}
	where $D=\rho^{-1}(\wt D\cap(\Phi(n)\setminus\Ro_n))$ and $\vfi(i,j)=\wt\vfi(\rho(i,j))$ (in this case, the elements $(2,1)$ and $(n,n-1)$ are not contained in $\wt D$, and the equations on the $(n-1)$th row are obtained from compatibility of system (\ref{syst:abkappa})).
	
	From the requirement that $\wt x\in\Supp{\psi}$ we also get that $M_{n,j}^{\wt D}(\wt x)=0$ when $2<j<n$ and $j\notin\Ro_{\wt D}$. It follows that $t_{n-1,j}=-\dfrac{x_{n,j}}{x_{n,n-1}}$ for $\vartheta<j<n-1$, and $t_{n-1,j}$ for $2<j<\vartheta$ can be expressed as a function of other $t_{n-1,j}$ as follows:
	\begin{equation}
		t_{n-1,j}=\eth-\frac{\sum\limits_{j<k\leq\vartheta}t_{n-1,k}M_{\wt\Ro_D}^{\wt X^k}(x)}{M_{\wt\Ro_D}^{\wt\Co_D}(x)}.
		\label{expr:h_n-1,j}
	\end{equation}
	Here, $\wt\Co_D=\{s\in\Co_D\mid s>j\}$, $\wt\Ro_D=\{s\in\Ro_D\mid s<n-j+1\}$, $\wt X^k$ is obtained from $\wt\Co_D$ by replacing $k$ by $j$, and $\eth$ is a certain constant. Solving system (\ref{syst:abkappa}) (and using expressions (\ref{expr:h_n-1,j})), we see that $\alpha_j(x)=0$ for $(i,j)\in D$.

	Substituting expression (\ref{expr:h_n-1,j}) into formula (\ref{formula:char}) and using previous arguments, we see that
	\begin{equation}
		\begin{split}
			&x_{j,1}-\sum\limits_{2<k<j, k\notin\Co_{D^+}}\frac{\text{sgn}(\Tu_{\wt X_{D^+}^{k}})x_{k,1}M_{\wt\Ro_{D^+}^k}^{\wt X^j_{D^+}}(x)}{M_{\wt\Ro_{D^+}^k}^{\wt\Co_{D^+}^k}(x)}=0\text{ for }(i,j)\in D^+,\\
			&\xi_{3,n}x_{n,n-1}=\xi_{1,n-1}x_{3,1}.
			\label{equations:3case_2part}
		\end{split}
	\end{equation}
	Here, $\wt\Co_{D^+}^{k}=\{s\in\Co_D\mid s>k\}$, $\wt\Ro_{D^+}^{k}=\{s\in\Ro_D\mid s<n-k+1\}$, $\wt X^k_{D^+}$ is obtained from $\wt\Co_{D^+}^{k}$ by replacing $j$ by $k$.
	
	Hence, in this case, $\Supp{\chi}$ is contained in the union of the subvarieties $P_D(\vfi)$ defined by equations (\ref{equations:3case_1part}) and (\ref{equations:3case_2part}). Again, such a subvariety $P_D(\vfi)$ is a union of conjugacy classes. More precisely,\break $P_D(\vfi)=\bigsqcup\limits_{\Do,\vfi}\Ko_{\Do}(\vfi)$, where the set $\Do\in\Phi(n)$ is obtained from the set $D$ by adding the elements $(3,1)$ and $(n,n-1)$, so that $\Do\supset D_3$. 
	
	The map~$\vfi$ is extended to a map $\Do\to\Fp_q^{\times}$ in an arbitrary way so that $\xi_{2,n-2}\vfi(3,1)=\xi_{3,n}\vfi(n,n-1)$. (This condition follows from the last equation in system (\ref{equations:3case_2part}) and requirements on minors in the definition of $\Ko_{\Do}(\vfi)$.) Thus, to conclude the proof in the case $x_{3,2}=0 $ and $x_{n,n-1}\neq 0$, it remains to check that $\chi(x)\neq0$ for all $x\in\Ko_{\Do}(\vfi)$. This will be done in the next subsection.
	
	At last, consider the case when $x_{3,2}\neq 0$ and $x_{n,n-1}\neq 0$. Using argumentation as in the previous case, we obtain
	
	\begin{equation*}
		\begin{split}
			&M^{\wt D}_{i,j}(\wt x)=M^{\wt D}_{i,j}(\wt x_{\wt D}(\vfi))\text{ for }(i,j)\in R(\wt D),\\
			&\beta_i^2(\wt x)=\alpha_j(\wt x)=0\text{ for }(i,j)\in \wt D\setminus((n-1,1)\cup(n,2)),\\
			&\xi_{2,n-2}\wt x_{3,2}=\xi_{3,n}\wt x_{n,n-2}.
		\end{split}  
	\end{equation*}
	But $\wt x_{i,j}=x_{i,j}$ when $i<n-1$ for all $t\in\Tu$. It follows that
	\begin{equation}
		\begin{split}
			&M_{i,j}^{D}(x)=M_{i,j}^{D}(x_D(\vfi))\text{ for }(i,j)\in R(D)\setminus(\Co_1\cup\Ro_n),\\
			&\beta_i^2(x)=0\text{ for }(i,j)\in D^+,
		\end{split}
		\label{equations:4case_1part}
	\end{equation}
	where $D=\rho^{-1}(\wt D\cap(\Phi(n)\setminus\Ro_n))$ and $\vfi(i,j)=\wt\vfi(\rho(i,j))$ (in this case, the elements $(2,1)$ and $(n,n-1)$ are contained in $\wt D$, and the equations on the $(n-1)$th row are obtained from the compatibility of system~(\ref{syst:abkappa})). Also we get that 
	\begin{equation*}
		\begin{split}
			&\alpha_j(x)-x_{n,n-1}x_{n-1,j}=-x_{n,n-1}\sum\limits_{k=\wp'}^{n-2}t_{n-1,k}x_{k,j}\text{ for }j\in D^+,\\
			&t_{n-1,n-2}=\cfrac{x_{3,2}\xi_{2,n-2}-x_{n,n-2}\xi_{3,n}}{x_{n,n-1}\xi_{3,n}}.
		\end{split}
	\end{equation*}
	Here $\wp'=\text{min}_{(b,s)\in D^+}b$. And it is easy to see from system (\ref{syst:abkappa}) that this system of equations (except the last one) is equivalent to the requirement that $\alpha_j(x)=0$ for $(i,j)\in D^+$. 
	
	Solving system (\ref{syst:abkappa}), we see that
	\begin{equation*}\begin{split}
			&t_{n-1,j}=\dfrac{\eth-\sum\limits_{\wp'<k<n-1,~k\notin\Ro_D}\text{sgn}(\tau_{X_k})t_{n-1,k}M^{\Co_D\setminus\{2\}}_{X_k}(x)}{M_{\Ro_D\setminus \{3\}}^{\Co_D\setminus \{2\}}(x)}\text{ for }j\in\Ro_D\setminus \{3\},\\
			&t_{n-1,3}=\dfrac{\eth-\sum\limits_{\wp'<k<n-1,~k\notin\Ro_D}\text{sgn}(\tau_{X_k})t_{n-1,k}M^{\Co_D}_{W_k}(x)}{M_{\Ro_D}^{\Co_D}(x)}-\dfrac{\sum\limits_{3<j<\wp'}t_{n-1,j}x_{j,2}}{x_{3,2}}.
		\end{split}    
	\end{equation*}
	Here the set $X_k$ is obtained from $\Ro_D\setminus\{3\}$ by replacing $j$ by $k$, $W_k$ is obtained from $\Ro_D\setminus \{3\}$ by adding $k$, and $\eth$ is a certain constant. Arguing similarly to the previous case, we deduce that, after substitutions of all these expressions into formula (\ref{formula:char}), the coefficients at the remaining variables should be zero. Thus, we obtain that
	\begin{equation*}
		\begin{split}
			&x_{2,1}=0,~-x_{i,1}x_{3,2}+x_{3,1}x_{i,2}=0,\text{ for }3<i<\wp',\\
			&x_{i,1}M_{\Ro_D\setminus \{3\}}^{\Co_D\setminus \{2\}}(x)-\sum\limits_{k\in\Ro_D}\text{sgn}(\tau_{X^k})x_{k,1}M_{X^k}^{\Co_D\setminus \{2\}}(x)-\frac{\text{sgn}(\tau_{W_i})x_{3,1}M_{W_i}^{\Co_D}}{x_{3,2}}=0,~\wp'<i<n-1,i\notin\Ro_D,
		\end{split}
	\end{equation*}
	
	These equations can be rewritten as follows:
	\begin{equation}
		M_{i,1}^{D}(x)=0,~1<i<n-2,i\notin\Ro_D.
		\label{equations:4case_2part}
	\end{equation}
	
	Hence, in this case, $\Supp{\chi}$ is contained in the union of the subvarieties $P_D(\vfi)$ defined by equations (\ref{equations:4case_1part}), (\ref{equations:4case_2part}) and $\alpha_j(x)=0$ for $(i,j)\in D^+$. Such a subvariety $P_D(\vfi)$ is a union of conjugacy classes. More precisely, $P_D(\vfi)=\bigsqcup\limits_{\Do_4,\vfi}\Ko_{\Do_4}(\vfi)\cup\bigsqcup\limits_{\Do_5,\vfi}\Ko_{\Do_5}(\vfi)$, where the set $\Do_4\in\Phi(n)$ is obtained from the set~$D$ by adding the element $(n,n-1)$ and a (probably, empty) subset of the set $\{(n-1,1),(n,3)\}$, so that $\Do_4\supset D_4$, while the set $\Do_5\in\Phi(n)$ is obtained from the set $D$ by adding the elements $(n,n-1)$, $(n-2,1)$, $(n-1,2)$ and, probably, $(n,3)$, so that $\Do_5\supset D_5$. 
	
	The map $\vfi$ is extended to a map $\Do_i\to\Fp_q^{\times}$, where $i=4,5$ in an arbitrary way so that\break $\xi_{3,n}^2\vfi(n-1,2)\vfi(n,n-1)^2-\xi_{2,n-2}\xi_{1,n-1}\vfi(n-2,1)\vfi(3,2)^2=0$ for $\Do_5$. (In fact, the last condition defines the value $\vfi(n-1,2)$.) Thus, to conclude the proof in the case $x_{3,2}\neq 0 $ and $x_{n,n-1}\neq 0$, it remains to check that $\chi(x)\neq0$ for all $x\in\Ko_{\Do_i}(\vfi)$ where $i=4,5$. This will be done in the next subsection.

	\sst{Proof of Theorem~\ref{theo:2_depth_value}}  We need to calculate the value of the character~$\chi$ corresponding to an arbitrary 2-regular orbit $\Omega_{D^2,\xi}$ on all conjugacy classes $\Ko_D(\vfi)$. Of course, we will calculate the value of $\xi$ on the element $x_D(\vfi)\in\Ko_D(\vfi)$. We will consider all cases $D\supset D_i$, $1\leq i\leq 5$, subsequently.
	
	First, let the subset $D$ contain $D_1$ and $\vfi\colon D\to\Fp_q^{\times}$ be an arbitrary map. Denote $x_D(\vfi)$ by $x=(x_{i,j})$ and $x^t_{\wt U^{\kappa}}$ by $\wt x$ for simplicity (recall the notion $x^t=txt^{-1}$ for $t\in\Tu$ from page~\pageref{page:x_t}). According to Proposition~\ref{prop:chi_2_regular_Mackey} and the standard formulas for induced characters, one has
	\begin{equation}\label{form.char}
		\chi(x)=\sum\limits_{\{t\in\Tu,~txt^{-1}\in AB^{\kappa}\}}\theta(\wt x_{n-1,1}\xi_{1,n-1})\psi(\wt x).
	\end{equation}
	Since $t\in\Tu$ satisfies $txt^{-1}\in AB^{\kappa}$, we have $t_{n-1,j}=0$ for $(j,n-j+1)\in\overline{D}$ (see formula~(\ref{formula:x_t}) for $x^t$). Using the formula for a character of depth 1 (see Theorem~\ref{theo:0_1_depth_character}), we can rewrite formula (\ref{form.char}) as follows:
	\begin{equation*}
		\chi(x)=\theta(x_{n-1,1}\xi_{1,n-1})\sum\limits_{\{t\in\Tu,~t_{n-1,j}=0\text{ for } (j,n-j+1)\in \overline{D}\}}\psi(\wt x)=q^{\wt m_D}\prod\limits_{(i,j)\in D}\theta(x_{i,j}\xi_{j,i}),
	\end{equation*}
	where
	\begin{equation*}
		\wt m_D=\begin{cases}|R(D)\cap\Phi_{\mathrm{2reg}}|-1+|D^{+}|-|\overline{D}|,&\text{if } \overline{D}=D^+\cup\{(n-2,2)\},\\ |R(D)\cap\Phi_{\mathrm{2reg}}|-2,&\text{otherwise}.
		\end{cases}
	\end{equation*}
	Now, formula~(\ref{formula:m_D_depth_2}) shows that $\wt m_D=m_D$, as required.
	
	Second, let the subset $D$ contain $D_2$. This case is similar to the previous one. We have $t_{n-1,3}=0$ and $t_{n-1,j}=0$ for $(j,n-j+1)\in D^{+}$, so formula (\ref{form.char}) can be rewritten as follows:
	\begin{equation*}
		\chi(x)=\theta(x_{n-1,1}\xi_{1,n-1})\sum\limits_{\{t\in\Tu,~t_{n-1,3}=0\text{, } t_{n-1,j}=0\text{ for } (j,n-j+1)\in D^{+}\}}\psi(\wt x)=q^{\wt m_D}\prod\limits_{(i,j)\in D}\theta(x_{i,j}\xi_{j,i}).
	\end{equation*}
	Here $$\wt m_D=|R(D)\cap\Phi^{'}|-(n-4-|D^{+}|)+2n-7-|D^{+}|=|R(D)\cap\Phi^{'}|+n-3.$$
	By formula~(\ref{formula:m_D_depth_2}), this number equals $m_D$, which completes the proof in this case.
	
	Next, in the case when $D$ contains $D_3$, formula (\ref{form.char}) gives us that
	\begin{equation*}
		\chi(x)=\sum\limits_{\{t\in M,~txt^{-1}\in AB^{\kappa}\}}\theta(-x_{3,1}\xi_{1,n-1}t_{n-1,3})\psi(\wt x).
	\end{equation*}
	We know that $t_{n-1,j}=0$ for $\vartheta<j<n-1$ because $\wt\chi(\wt x)\neq 0$. The determinants $M_{i,j}^{D}(txt^{-1})$ are zero for $(n,j)\in R(D)$, $3<j<\vartheta$, by the same reason. It follows that $t_{n-1,j}=0$ for $(n,j)\in R(D)$, $3<j<\vartheta$, thus,
	\begin{equation*}
		\begin{split}
			\chi(x)&=q^{m'_D}\prod\limits_{(i,j)\in D}\theta(x_{i,j}\xi_{j,i})\quad
			\sum\limits_{t\in N}\theta(-x_{3,1}\xi_{1,n-1}t_{n-1,3})\theta\left(\xi_{3,n}\frac{(-1)^{|D^+|\text{ (mod 2)}}M_{n,3}^{D}(txt^{-1})}{M_{n-\ell+1,\ell}^{D}(x)}\right)\\
			&=q^{m'_D}\prod\limits_{(i,j)\in D}\theta(x_{i,j}\xi_{j,i})\quad
			\sum\limits_{t\in N}\theta(-x_{3,1}\xi_{1,n-1}t_{n-1,3})\theta(\xi_{3,n}x_{n,n-1}t_{n-1,3}+x_{n,3}\xi_{3,n})\\
			&=q^{\wt m_D}\prod\limits_{(i,j)\in D}\theta(x_{i,j}\xi_{j,i}).
		\end{split}
	\end{equation*}
	Here $m'_D$ is the constant defined by formula~(\ref{formula:m_D_0_1}) for subregular characters, $$N=\{t\in \Tu\mid t_{n-1,2}\text{, }t_{n-1,3}\neq 0\text{, } t_{n-1,j}\neq 0\text{ if }(n-j+1,j)\in D^{+}\},$$ and, clearly, $$\wt m_D=m'_D+|N|=|R(D)\cap\Phi^{''}|-1+n-4-|D^{+}|+|D^{''}|+2=|R(D)\cap\Phi^{''}|+n-3.$$ The latter number equals $m_D$, so the proof in this case is complete.

	Next, let $D$ contain $D_4$. From formula (\ref{form.char}) we see that
	\begin{equation*}
		\chi(x)=\sum\limits_{\{t\in M,~txt^{-1}\in AB^{\kappa}\}}\theta(x_{n-1,1}\xi_{1,n-1})\psi(\wt x).
	\end{equation*}
	The requirement that $x^t\in AB^{\kappa}$ gives us $t_{n-1,j}=0$ for $j\in \Ro_{\overline{D}}$, at the same time, the requirement that $\wt x\in\Supp{\psi}$ gives us $t_{n-1,n-2}=\cfrac{x_{3,2}\xi_{2,n-2}}{x_{n,n-1}\xi_{3,n}}$, so we can say that $\wt\gamma(\wt x)=\wt\gamma(x_{\wt D}(\wt\vfi))$ (here $\wt\gamma$ is the polynomial $\gamma$ from (\ref{formula:alpha_j_beta_i})). Hence we can rewrite formula (\ref{form.char}) as follows:
	\begin{equation*}\begin{split}
			\chi(x)&=q^{m'_D}\prod\limits_{(i,j)\in \overline{D}}\theta(x_{i,j}\xi_{j,i})\sum\limits_{t\in N}\theta(\xi_{1,n-1}x_{n-1,1})\theta(\vfi(n,3)\xi_{3,n}+\vfi(n-2,2)\xi_{2,n-2})\\
			&=q^{m'_D}\prod\limits_{(i,j)\in \overline{D}}\theta(x_{i,j}\xi_{j,i})\sum\limits_{t\in N}\theta(\xi_{1,n-1}x_{n-1,1}+x_{n,3}\xi_{3,n})\\
			&=q^{\wt m_D}\prod\limits_{(i,j)\in D}\theta(x_{i,j}\xi_{j,i}).
		\end{split}
	\end{equation*}
	
	Here $m'_D$ is the constant defined by formula~(\ref{formula:m_D_0_1}) for subregular characters, $\vfi(i,j)=\wt\vfi(\rho(i,j))$ and $$N=\{t\in\Tu\mid t_{n-1,j}=0\text{ for }j\in\Ro_{\overline{D}}\text{, }t_{n-1,n-2}=\cfrac{x_{3,2}\xi_{2,n-2}}{x_{n,n-1}\xi_{3,n}}\},$$ and it is easy to check that $$\wt m_D=m'_D+|N|=|R(D)\cap\Phi^{'}|+n-2+|D^{+}|-n+4+n-3-|D^{+}|-2=|R(D)\cap\Phi^{''}|+n-3.$$ The latter number equals $m_D$, so the proof in this case is complete.

	Finally, consider the case when $D$ contains $D_5$. From formula (\ref{form.char}) we see that
	\begin{equation*}
		\chi(x)=\sum\limits_{\{t\in M,~txt^{-1}\in AB^{\kappa}\}}\theta(-x_{n-2,1}\xi_{1,n-1}t_{n-1,n-2})\psi(\wt x).
	\end{equation*}
	The requirement that $x^t\in AB^{\kappa}$ gives us $t_{n-1,j}=0$ for $j\in \Ro_{D^+}$ and $t_{n-1,3}=\cfrac{x_{n-1,2}}{x_{3,2}}$, at the same time, the requirement that $\wt x\in\Supp{\psi}$ gives us $t_{n-1,n-2}=\cfrac{x_{3,2}\xi_{2,n-2}}{x_{n,n-1}\xi_{3,n}}$, so we can say that $\wt\gamma(\wt x)=\wt\gamma(x_{\wt D}(\wt\vfi))$ (here, $\wt\gamma$ is the polynomial $\gamma$ from (\ref{formula:alpha_j_beta_i})). Hence we can rewrite formula (\ref{form.char}) as follows:
	\begin{equation*}\begin{split}
			\chi(x)&=q^{m'_D}\prod\limits_{(i,j)\in \overline{D}}\theta(x_{i,j}\xi_{j,i})\sum\limits_{t\in N}\theta(-\xi_{1,n-1}\frac{x_{n-2,1}x_{3,2}\xi_{2,n-2}}{x_{n,n-1}\xi_{3,n}})\theta(\vfi(n,3)\xi_{3,n}+\vfi(n-2,2)\xi_{2,n-2})\\
			&=q^{m'_D}\prod\limits_{(i,j)\in \overline{D}}\theta(x_{i,j}\xi_{j,i})\sum\limits_{t\in N}\theta(\xi_{3,n}x_{n,3}+\frac{x_{n,n-1}x_{n-1,2}\xi_{3,n}}{x_{3,2}})\\
			&=q^{\wt m_D}\prod\limits_{(i,j)\in D}\theta(x_{i,j}\xi_{j,i}).
		\end{split}
	\end{equation*}
	
	Here $m'_D$ is the constant defined by formula~(\ref{formula:m_D_0_1}) for subregular characters, $\vfi(i,j)=\wt\vfi(\rho(i,j))$, $$N=\left\{t\in\Tu\mid t_{n-1,j}=0\text{ for }j\in\Ro_{D^+}\text{, }t_{n-1,3}=\cfrac{x_{n-1,2}}{x_{3,2}}\text{, }t_{n-1,n-2}=\cfrac{x_{3,2}\xi_{2,n-2}}{x_{n,n-1}\xi_{3,n}}\right\},$$ and it is easy to check that $$\wt m_D=m'_D+|N|=|R(D)\cap\Phi^{'}|+n-2+|D^{+}|-n+4+n-3-|D^{+}|-2=|R(D)\cap\Phi^{''}|+n-3.$$ The latter number equals $m_D$, so the proof in this case is complete.

	\bigskip\textsc{Mikhail V. Ignatev: National Research University Higher School of Economics, Pokrovsky Boulevard 11, 109028, Moscow, Russia}
	
	\emph{E-mail address}: \texttt{mihail.ignatev@gmail.com}
	
	\bigskip\textsc{Mikhail S. Venchakov: St. Petersburg University, Universitetskaya Embankment 7--9, 199034, St. Petersburg, Russia}
	
	\emph{E-mail address}: \texttt{mihail.venchakov@gmail.com}
	

\begin{thebibliography}{XXXX}\addcontentsline{toc}{section}{References}
		
		\bibitem[An95]{Andre95} C.A.M. Andr\`e. Basic sums of coadjoint orbits of the unitriangular group. J.~Algebra \textbf{176} (1995), 959--1000.
		
		\bibitem[An01]{Andre01} C.A.M. Andr\`e. The basic character table of the unitriangular group. J.~Algebra \textbf{241} (2001), 437--471.
		
		
		
		
		
		
		
		
		
		
		
		\bibitem[Ig09]{Ignatev09} M.V. Ignatev. Subregular characters of the unitriangular group over a finite field. J. Math. Sci. \textbf{156} (2009), no. 2, 276--291, arXiv: \texttt{math.RT/0801.3079}.
		
		
		\bibitem[IPa09]{IgnatevPanov09} M.V. Ignatev, A.N. Panov. Coadjoint orbits of the group $\mathrm{UT}(7, K)$. J. Math. Sci. \textbf{156} (2009), no. 2, 292--312; arXiv: \texttt{math.RT/0603649}.
		
		
		
		
		\bibitem[Ka77]{Kazhdan77} D. Kazhdan. Proof of Springer's hypothesis.
		Israel J. Math. \textbf{28} (1977), 272--286
		
		\bibitem[Ki62]{Kirillov62} A.A. Kirillov. Unitary representations of nilpotent Lie groups. Russian Math. Surveys \textbf{17}~(1962), 53--110.
		
		\bibitem[Ki04]{Kirillov04} A.A. Kirillov. Lectures on the orbit method. Grad. Stud. in Math. \textbf{64}, AMS, 2004.
		
		
		\bibitem[Kr00]{Kraft00} H. Kraft. Geometric methods in invariant theory. IO NFMI, 2000.
		
		
		\bibitem[Le74]{Lehrer74} G.I. Lehrer. Discrete series and the unipotent subgroup. Compositio Math., \textbf{28} (1974),\break fasc. 1, 9--19.
		
		
		
		
		
		\bibitem[Pa08]{Panov08} A.N. Panov. Involutions in $S_n$ and associated coadjoint orbits. J. Math. Sci. \textbf{151} (2008), 3018--3031.
		
		\bibitem[Pa09]{Panov09} A.N. Panov. On the index of certain nilpotent Lie algebras. J. Math. Sci. \textbf{161} (2009), no. 1, 122--129.
		
		\bibitem[Ve70]{Vergne} M. Vergne. Construction de sous-alg\`{e}bres subordonn\'{e}es \`{a} un \'{e}l\'{e}ment du~dual d'une alg\`{e}bre de Lie r\'{e}soluble. C. R. Acad. Sci. Paris Ser. A--B \textbf{270} (1970), A173--A175.
		
		
	\end{thebibliography}
\end{document}